\documentclass[a4paper, 11pt, reqno]{amsart}

\usepackage{amsmath, latexsym, amsfonts, amssymb, amsthm, amscd,
  stmaryrd}%
\usepackage{graphicx, color}%
\usepackage{natbib}%
\usepackage{setspace}

\definecolor{darkblue}{rgb}{0.0,0.0,0.7}

\usepackage[dvips,%
bookmarks = true,%
colorlinks = true,%
linkcolor = darkblue,%
citecolor = darkblue,%
urlcolor = darkblue, %
]{hyperref}

\numberwithin{equation}{section}

\newtheorem{theorem}{Theorem} 
\newtheorem{lemma}{Lemma}
\newtheorem{proposition}{Proposition}
\newtheorem{corollary}{Corollary}

\theoremstyle{remark}
\newtheorem*{remark}{Remark}

\theoremstyle{definition}
\newtheorem{definition}{Definition}

\theoremstyle{remark}%
\newtheorem*{example}{Example}

\newcounter{assumption}
\newtheorem{assumption}[assumption]{Assumption}

\renewcommand{\leq}{\leqslant}
\renewcommand{\geq}{\geqslant}

\newcommand{\Span}{\text{Span}}

\newcounter{hyp}

\newcommand{\mc}{\mathcal}
\newcommand{\mb}{\mathbf}
\newcommand{\mbb}{\mathbb}
\newcommand{\mrm}{\mathrm}
\newcommand{\mf}{\mathfrak}
\newcommand{\msf}{\mathsf}

\newcommand{\tup}[1]{\textup{#1}}

\newcommand{\wh}{\widehat}

\newcommand{\von}{\varepsilon}

\newcommand{\tta}{\theta}

\newcommand{\Lba}{\Lambda}

\newcommand{\lba}{\lambda}

\newcommand{\raro}{\rightarrow}
\newcommand{\sumin}{\sum_{i=1}^{n}}
\newcommand{\setR}{\mathbb{R}}

\newcommand{\ind}[1]{\mathbf{1}_{#1}}

\DeclareMathOperator*{\argmin}{argmin}
\DeclareMathOperator*{\argmax}{argmax}

\newcommand{\prodsca}[2]{\langle #1 \, ,\, #2 \rangle}

\newcommand{\norm}[1]{\| #1 \|}

\newcommand{\norminfty}[1]{\| #1 \|_{\infty}}

\newcommand{\ppint}[1]{\lfloor #1 \rfloor} 

\DeclareMathOperator*{\trace}{Tr}%

\newcommand{\diag}{\text{diag}}

\newlength{\figurelength}%

\title{Uniform estimation of a signal based on inhomogeneous data}

\author{St{\'e}phane Ga{\"\i}ffas}

\address{Modal'X, Universit\'e Paris X -- Nanterre, b\^atiment G, 200
  avenue de la R\'epublique, 92000 Nanterre }

\email{\href{mailto:stephane.gaiffas@u-paris10.fr}{
    \texttt{stephane.gaiffas@u-paris10.fr}}}



\date{\today}

\dedicatory{ Modal'X, Universit\'e Paris X -- Nanterre \\
  email\tup: \rm{ \href{mailto:stephane.gaiffas@u-paris10.fr}{
      \texttt{stephane.gaiffas@u-paris10.fr} } } }

\begin{document}

\begin{abstract}
  We want to reconstruct a signal based on inhomogeneous data (the
  amount of data can vary strongly), using the model of regression
  with a random design.  Our aim is to understand the consequences of
  inhomogeneity on the accuracy of estimation within the minimax
  framework. Using the uniform metric weighted by a
  spatially-dependent rate as a benchmark for an estimator accuracy,
  we are able to capture the deformation of the usual minimax rate in
  situations with local lacks of data (modelled by a design density
  with vanishing points). In particular, we construct an estimator
  both design and smoothness adaptive, and a new criterion is
  developed to prove the optimality of these deformed rates.
\end{abstract}

\maketitle

\section{Introduction}
\label{sec:intro}

\subsection*{Motivations}

A problem particularly prominent in statistical literature is the
adaptive reconstruction of a function based on irregularly sampled
noisy data. In several practical situations, the statistician cannot
obtain ``nice'' regularly sampled observations, because of various
constraints linked with the source of the data, or the way the data is
obtained. For instance, in signal or image processing, the irregular
sampling can be due to the process of motion or disparity compensation
(used in advanced video processing), while in topography, measurement
constraints are linked with the properties of the
ground. See~\cite{grochenig94} for a survey on irregular sampling,
\cite{almansa_rouge_jaffard03}, \cite{konrad00} for applications
concerning respectively satellite image and stereo imaging, and
\cite{jansen04} for examples of geographical constraints.

Such constraints can result in potentially strong local lacks of
data. Consequently, the accuracy of a procedure based on such data can
become locally very poor. The aim of the paper is to study from a
theoretical point of view the consequences of data
\emph{inhomogeneity} on the reconstruction of a univariate
signal. Natural questions arise: how does the inhomogeneity impact on
the accuracy of estimation?  What does the optimal convergence rate
become in such situations? Can the rate vary strongly from place to
place, and how?

\subsection*{The model}

The widest spread way to model such observations is as follows. We
model the available data $[(X_i, Y_i); 1 \leq i \leq n]$ by
\begin{equation}
  \label{eq:model}
  Y_i = f(X_i) + \sigma \xi_i,
\end{equation}
where $\xi_i$ are i.i.d. Gaussian standard and independent of the
$X_i$'s and $\sigma > 0$ is the noise level. The design variables
$X_i$ are i.i.d. with unknown density $\mu$ on $[0, 1]$. The more the
density $\mu$ is ``far'' from the uniform law, the more the data drawn
from~\eqref{eq:model} is inhomogeneous. A simple way to include
situations with local lacks of data within the model~\eqref{eq:model}
is to allow the density $\mu$ to be arbitrarily small at some points,
and to vanish. This kind of behaviour is not commonly used in
literature, since most papers assume $\mu$ to be uniformly bounded
away from zero. We give references handling this kind of design below.

In practice, we don't know $\mu$, since it requires to know in a
precise way the constraints making the observation irregularly
sampled, neither do we know the smoothness of $f$. Therefore, a
convenient procedure shall adapt both to the design and to the
smoothness of $f$. Such a procedure (that is proved to be optimal) is
constructed here.

\subsection*{Methodology}

We want to reconstruct $f$ globally, with sup norm loss. The reason
for choosing this metric is that it is exacting: roughly, it forces an
estimator to behave well at every point simultaneously. This property
is convenient here, since it allows to capture in a very simple way
the consequences of inhomogeneity directly on the convergence rate.

In what follows, $a_n \lesssim b_n$ means $a_n \leq C b_n$ for any
$n$, where $C > 0$. We say that a sequence of curves $v_n(\cdot) \geq
0$ is an upper bound over some class $F$ if there is an estimator $\wh
f_n$ such that
\begin{equation}
  \label{eq:new_UB}
  \sup_{f \in F} \mb E_{f\mu} \big[ w\big( \sup_{x \in [0, 1]}
  v_n(x)^{-1} |\wh f_n(x) - f(x) | \big) \big] \lesssim 1
\end{equation}
as $n \raro +\infty$, where $\mb E_{f\mu}$ denotes the expectation
with respect to the joint law $\mb P_{f\mu}$ of the $[(X_i, Y_i); 1
\leq i \leq n]$, and where $w(\cdot)$ is a loss function, that is a
non-negative and non-decreasing function such that $w(0) = 0$ and
$w(x) \leq A (1 + |x|^b)$ for some $A, b > 0$.

\subsection*{Literature}

Pointwise estimation at a point where the design vanishes is studied
in~\cite{hall_et_al97}, with the use of a local linear procedure. This
design behaviour is given as an example in~\cite{guerre99}, where a
more general setting for the design is considered, with a Lipschitz
regression function. In~\cite{gaiffas04a}, pointwise minimax rates
over H\"older classes are computed for several design behaviours, and
an adaptive estimator for pointwise risk is constructed
in~\cite{gaiffas05}. In these papers, it appears that, depending on
the design behaviour at the estimation point, the range of minimax
rates is very wide: from very slow (logarithmic) rates to very fast
quasi-parametric rates.

Many adaptive techniques have been developed in literature for
handling irregularly sampled data. Among wavelet methods, see
\cite{hall_et_al97} for interpolation; \cite{antoniadis_et_al97},
\cite{antoniadis_et_al98}, \cite{cai_brown98},
\cite{hall_park_turlach98}, \cite{wong_zheng02} for tranformation and
binning; \cite{antoniadis_fan01} for a penalization approach;
\cite{delouille_et_al01} and \cite{delouille_et_al04} for the
construction of design-adapted wavelet via lifting;
\cite{pensky_wiens01} for projection-based techniques and
\cite{kerk_picard_warped_04} for warped wavelets. For model selection,
see \cite{baraud02}. See also the PhD manuscripts from
\cite{voichitaphd} and \cite{delouille_phd}.

\section{Results}

To measure the smoothness of $f$, we consider the standard H\"older
class $H(s, L)$ where $s, L > 0$, defined as the set of all the
functions $f : [0, 1] \raro \setR$ such that
\begin{equation*}
  | f^{(\ppint{s})}(x) - f^{(\ppint{s})}(y) | \leq L |x - y|^{s -
    \ppint{s}}, \quad \forall x, y \in [0, 1],
\end{equation*}
where $\ppint{s}$ is the largest integer smaller than $s$. Minimax
theory over such classes is standard: we know from~\cite{stone82} that
within the model~\eqref{eq:model}, the minimax rate is equal to $(\log
n / n)^{s / (2s + 1)}$ over such classes, when $\mu$ is continuous and
uniformly bounded away from zero. If $Q > 0$, we define $H^Q(s, L) :=
H(s, L) \cap \{ f \;|\; \norminfty{f} \leq Q \}$ (the constant $Q$
needs not to be known).

We use the notation $\mu(I) := \int_I \mu(t) dt$. If $F = H(s, L)$ is
fixed, we consider the sequence of positive curves $h_n(\cdot) =
h_n(\cdot; F, \mu)$ satisfying
\begin{equation}
  \label{eq:r_n_def}
  L h_n(x)^s = \sigma \Big( \frac{\log n}{n \mu([x-h, x+h])}
  \Big)^{1/2}
\end{equation}
for any $x \in [0, 1]$, and we define
\begin{equation*}
  r_n(x; F, \mu) := L h_n(x; F,\mu)^s.
\end{equation*}
Since $h \mapsto h^{2s} \mu([x-h, x+h])$ is increasing for any $x$,
these curves are well-defined (for $n$ large enough) and unique. In
Theorem~\ref{thm:upper_bound} below, we show that $r_n(\cdot)$ is an
upper bound over H\"older classes, and the optimality of this rate is
proved in Theorem~\ref{thm:lower_bound}.
\begin{example}
  \newlength{\ecartt}%
  \setlength{\ecartt}{1cm}%
  When $s = 1$, $\sigma = L = 1$ and $\mu(x) = 4|x - 1/2|
  \ind{[0,1]}(x)$, solving~\eqref{eq:r_n_def} leads to
  \begin{equation*}
    r_n(x) = (\log n / n)^{\,\alpha_n(x)},
  \end{equation*}
  where the exponent $\alpha_n(\cdot)$ is given by
  \begin{equation*}
    \alpha_n(x) =
    \begin{cases}
      \, \frac{1}{3} \big(1 - \frac{\log( 1 - 2x )}{ \log (\log n /
        n)} \big) \hspace{\ecartt} \text{ when } x \in \big[0,
      \frac{1}{2} -
      ( \frac{ \log n}{2 n} )^{1/4} \big], \\[0.2cm]
      \, \frac{ \log \big( ( (x - 1/2)^4 + 4 \log n / n )^{1/2} -
        (x - 1/2)^2 \big) - \log 2 }{2 \log (\log n / n)} \\
      \phantom{\, \frac{1}{3} \big(1 - \frac{\log( 1 - 2x )}{ \log
          (\log n / n)} \big)} \hspace{\ecartt} \text{ when } x \in
      \big[ \frac{1}{2} - ( \frac{\log n}{2 n} )^{1/4}, \frac{1}{2} +
      ( \frac{\log n}{2 n} )^{1/4} \big], \\[0.2cm]
      \, \frac{1}{3} \big(1 - \frac{\log( 2 x - 1 )}{ \log(\log n /n)}
      \big) \hspace{\ecartt} \text{ when } x \in \big[\frac{1}{2} + (
      \frac{ \log n}{2 n} )^{1/4}, 1 \big].
    \end{cases}
  \end{equation*}
  Within this example, $r_n(\cdot)$ switches from one ``regime'' to
  another. Indeed, in this example there is a lack of data in the
  middle of the unit interval. The consequence is that $r_n(1/2) =
  (\log n / n)^{1/4}$ is slower than the rate at the boundaries
  $r_n(0) = r_n(1) = (\log n / n)^{1 / 3}$, which comes from the
  standard minimax rate $(\log n / n)^{s / (2s + 1)}$ with $s = 1$. We
  show the shape of this deformed rate for several sample sizes in
  Figure~\ref{fig:rate_example}.%
  \setlength{\figurelength}{3in}%
  \begin{figure}[htbp]
    \centering
    \includegraphics[width=\figurelength]{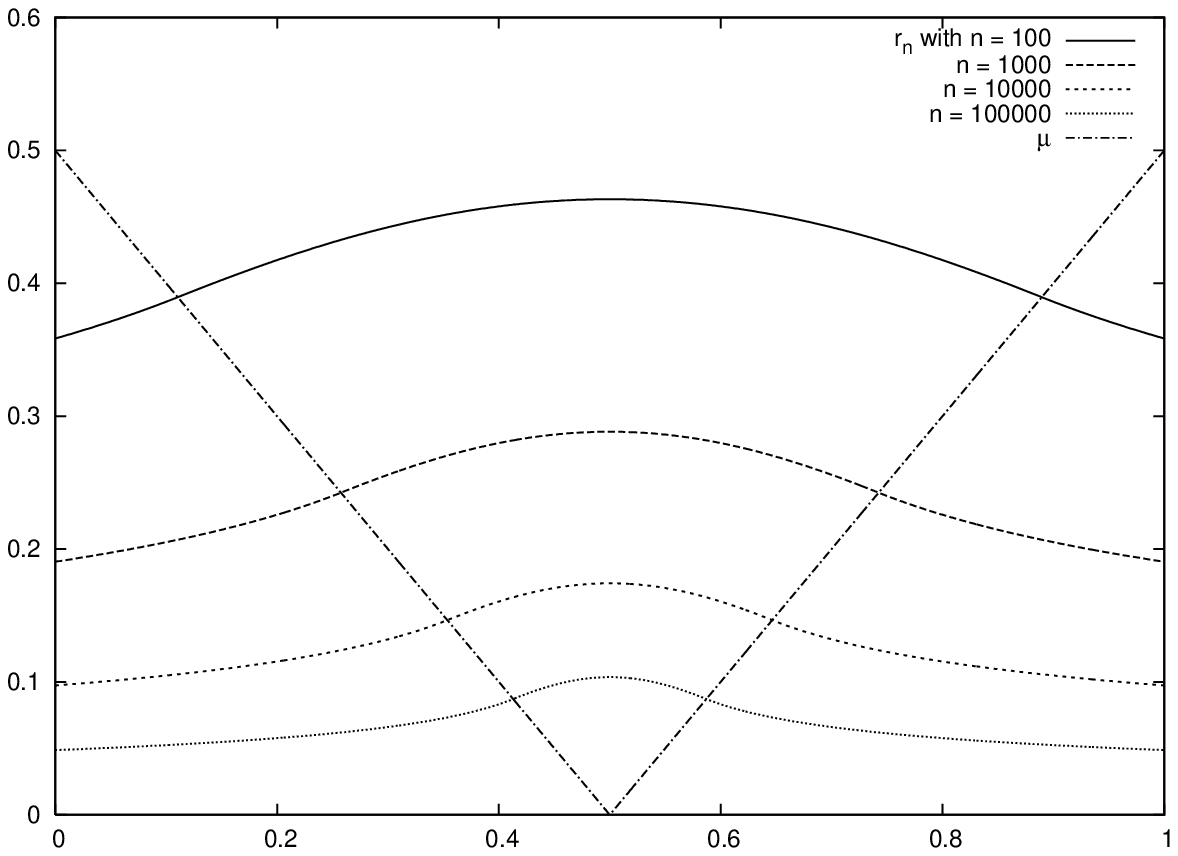}%
    \includegraphics[width=\figurelength]{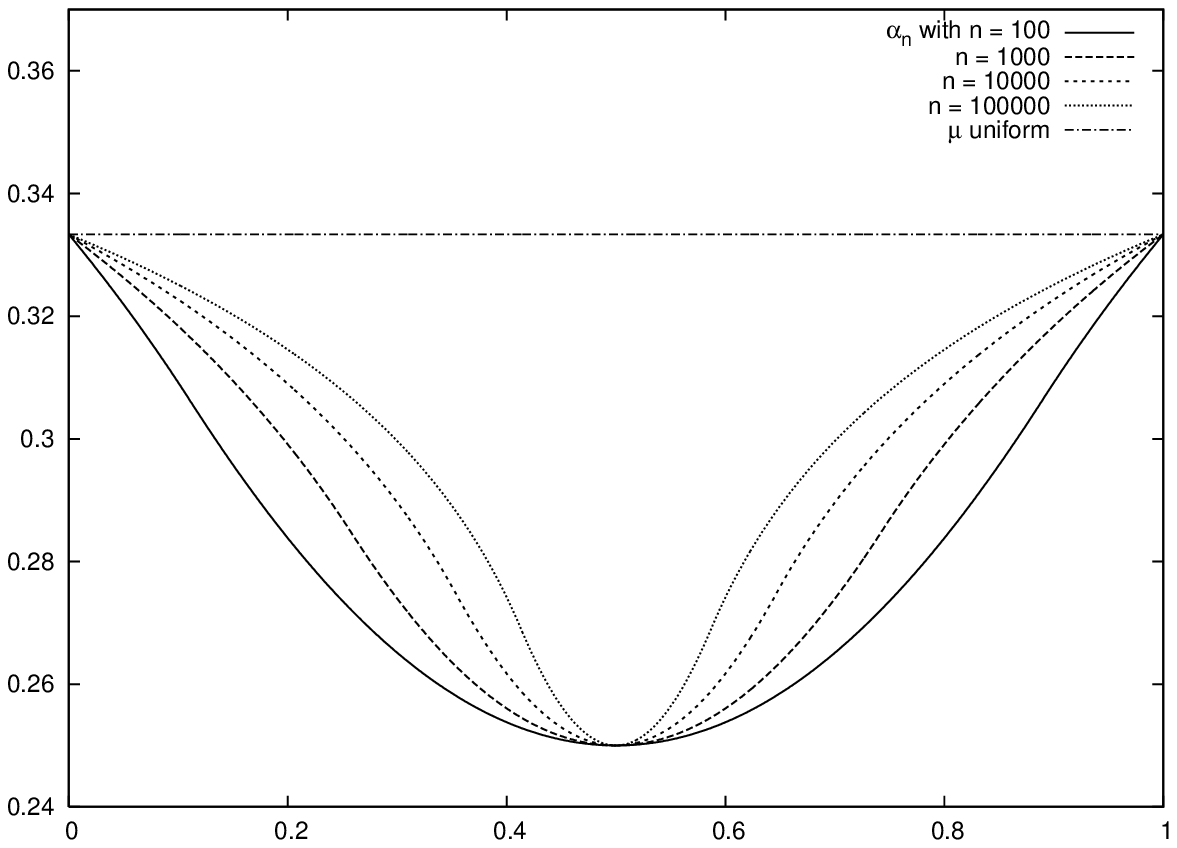}%
    \caption{$r_n(\cdot)$ and $\alpha_n(\cdot)$ for several sample
      sizes}
    \label{fig:rate_example}
  \end{figure}
\end{example}

\subsection*{Upper bound}

In this section, we show that the spatially-dependent rate
$r_n(\cdot)$ defined by~\eqref{eq:r_n_def} is an upper bound in the
sense of~\eqref{eq:new_UB} over H\"older classes. The estimator used
in this upper bound is both smoothness and design adaptive (it does
not depend on the design density within its construction). This
estimator is constructed in Section~\ref{sec:estimator} below. Let $R$
be a fixed natural integer.

\setcounter{assumption}{3}
\begin{assumption}
  \label{ass:design}
  We assume that $\mu$ is continuous, and that whether $\mu(x) > 0$
  for any $x$, or $\mu(x) = 0$ for a finite number of $x$. Moreover,
  for any $x$ such that $\mu(x) = 0$ we assume that $\mu(y) = |y -
  x|^{\beta(x)}$ for any $y$ in a neighbourhood of $x$ (where
  $\beta(x) \geq 0$).
\end{assumption}

\begin{theorem}
  \label{thm:upper_bound}
  Let $s \in (0, R+1]$ and assumption~\ref{ass:design} holds. The
  estimator $\wh f_n$ defined by~\eqref{eq:estimator} satisfies
  \begin{equation}
    \label{eq:upper_bound}
    \sup_{f \in F} \mb E_{f\mu} \big[ w( \sup_{x \in [0, 1]} r_n(x)^{-1}
    |\wh f_n(x) - f(x)| ) \big] \lesssim 1
  \end{equation}
  as $n \raro +\infty$ for any $F = H^Q(s, L)$, where $r_n(\cdot) =
  r_n(\cdot; F,\mu)$ is given by~\eqref{eq:r_n_def}.
\end{theorem}

This theorem assesses the adaptive estimator constructed in
Section~\ref{sec:estimator} below. The estimator $\wh f_n$ is based on
a precise estimation of the scaling coefficients (within a
multiresolution analysis) of $f$. This method relies on a Lepski-type
method (see for instance \cite{lepski_mammen_spok_97}) that we adapt
for random designs.

\begin{remark}
  Within Theorem~\ref{thm:upper_bound}, there are mainly two
  situations.
  \begin{itemize}
  \item $\mu(x) > 0$ for any $x$:~we have $r_n(x) \asymp (\log n /
    n)^{s / (2s + 1)}$ for any $x$, where $a_n \asymp b_n$ means $a_n
    \lesssim b_n$ and $b_n \lesssim a_n$. Hence, we find back the
    standard minimax rate in this situation. Note that this result is
    new since adaptive estimators over H\"older balls in regression
    with random design were not previously constructed.
  \item $\mu(x) = 0$ for one or several $x$:~the rate $r_n(\cdot)$ can
    vary strongly from place to place, depending on the behaviour of
    $\mu$. Indeed, the rate changes \emph{in order} from one point to
    another, see the example above.
  \end{itemize}
\end{remark}

\begin{remark}
  Implicitly, we assumed in Theorem~\ref{thm:upper_bound} that $s \in
  (0, R+1]$, where $R$ is a tuning parameter of the procedure. Indeed,
  in the minimax framework considered here, the fact of knowing an
  upper bound for $s$ is usual in the study of adaptive methods, and
  somehow, unavoidable. For instance, when considering adaptive
  wavelet methods, the ``maximum smoothness'' corresponds to the
  number of moments of the mother wavelet.
\end{remark}

\subsection*{Optimality of~$r_n(\cdot)$}

We have seen that the rate $r_n(\cdot)$ defined by~\eqref{eq:r_n_def}
is an upper bound over H\"older classes, see
Theorem~\ref{thm:upper_bound}. In Theorem~\ref{thm:lower_bound} below,
we prove that this rate is indeed optimal. In order to show that
$r_n(\cdot)$ is optimal in the minimax sense over some class $F$, the
classical criterion consists in showing that
\begin{equation}
  \label{eq:classical_lb}
  \inf_{\wh f_n} \sup_{f \in F} \mb E_{f\mu} \big[ w(\sup_{x \in [0,
    1]} r_n(x)^{-1} |\wh f_n(x) - f(x)|) \big] \gtrsim 1,
\end{equation}
where the infimum is taken among all estimators based on the
observations~\eqref{eq:model}. However, this criterion does not
exclude the existence of another normalisation $\rho_n(\cdot)$ that
can improve $r_n(\cdot)$ in some regions of $[0,
1]$. Indeed,~\eqref{eq:classical_lb} roughly consists in a minoration
of the uniform risk over the whole unit interval and then, only over
some particular points. Therefore, we need a new criterion that
strengthens the usual minimax one to prove the optimality of
$r_n(\cdot)$. The idea is simple: we localize~\eqref{eq:classical_lb}
by replacing the supremum over $[0, 1]$ by a supremum over any (small)
inverval $I_n \subset [0, 1]$, that is
\begin{equation}
  \label{eq:new_criterion}
  \inf_{\wh f_n} \sup_{f \in F} \mb E_{f\mu} \big[ w(\sup_{x \in I_n}
  r_n(x)^{-1} |\wh f_n(x) - f(x)|) \big] \gtrsim 1, \quad \forall I_n.
\end{equation}
It is noteworthy that in~\eqref{eq:new_criterion}, the length of the
intervals cannot be arbitrarily small. Actually, if an interval $I_n$
has a length smaller than a given limit,~\eqref{eq:new_criterion} does
not hold anymore. Indeed, beyond this limit, we can improve
$r_n(\cdot)$ for the risk localized over~$I_n$: we can construct an
estimator $\wh f_n$ such that
\begin{equation}
  \label{eq:improve_r_n}
  \sup_{f \in F} \mb E_{f\mu} \big[ w(\sup_{x \in I_n} r_n(x)^{-1}
  |\wh f_n(x) - f(x)|) \big] = o(1),
\end{equation}
see Proposition~\ref{prop:suroptimal} below. The phenomenon described
in this section, which concerns the uniform risk, is linked with the
results from~\cite{cai_low05} for shrunk $\mbb L^2$ risks. In what
follows, $|I|$ stands for the length of an interval~$I$.

\begin{theorem}
  \label{thm:lower_bound}
  Suppose that
  \begin{equation}
    \label{eq:assump_low_bound}
    \mu(I) \gtrsim |I|^{\beta + 1}
  \end{equation}
  uniformly for any interval $I \subset [0, 1]$, where $\beta \geq 0$
  and let $F = H(s, L)$. Then, for any interval $I_n \subset [0, 1]$
  such that
  \begin{equation}
    \label{eq:I_n_length}
    |I_n| \sim n^{-\alpha}
  \end{equation}
  with $\alpha \in (0, (1 + 2s + \beta)^{-1})$, we have
  \begin{equation}
    \label{eq:lower_bound}
    \inf_{\wh f_n} \sup_{f \in F} \mb E_{f\mu} \big[ w \big( \sup_{x
      \in I_n} r_n(x)^{-1} |\wh f_n(x) - f(x)| \big) \big] \gtrsim 1
  \end{equation}
  as $n \raro +\infty$, where $r_n(\cdot) = r_n(\cdot\;; F, \mu)$ is
  given by~\eqref{eq:r_n_def}.
\end{theorem}

\begin{corollary}
  \label{cor:lower_bound}
  If $v_n(\cdot)$ is an upper bound over $F = H(s, L)$ in the sense
  of~\eqref{eq:new_UB}, we have
  \begin{equation*}
    \sup_{x \in I_n} v_n(x) / r_n(x) \gtrsim 1
  \end{equation*}
  for any interval $I_n$ as in Theorem~\ref{thm:lower_bound}. Hence,
  $r_n(\cdot)$ cannot be improved uniformly over an interval with
  length $n^{\von -1 / (1 + 2s + \beta)}$, for any arbitrarily small
  $\von > 0$.
\end{corollary}

\begin{proposition}
  \label{prop:suroptimal}
  Let $F = H(s, L)$ and $\ell_n$ be a positive sequence satisfying
  \begin{equation*}
    \log \ell_n = o(\log n).
  \end{equation*}
  $a)$~Let $\mu$ be such that $0 < \mu(x) < +\infty$ for any $x \in
  [0, 1]$. Note that in this case, $r_n(x) \asymp (\log n / n)^{s /
    (2s + 1)}$ for any $x \in [0, 1]$ and
  that~\eqref{eq:assump_low_bound} holds with $\beta = 0$. If $I_n$ is
  an interval satisfying
  \begin{equation*}
    |I_n| \sim (\ell_n / n)^{1 / (1 + 2s)},
  \end{equation*}
  we can contruct an estimator $\wh f_n$ such that
  \begin{equation*}
    \sup_{f \in F} \mb E_{f\mu} \Big[ w \Big( \Big( \frac{n}{\log n}
    \Big)^{s / (2s + 1)} \sup_{x \in I_n} |\wh f_n(x) - f(x)| \Big)
    \Big] = o(1).
  \end{equation*}
  $b)$~Let $\mu(x_0) = 0$ for some $x_0 \in [0, 1]$ and $\mu([x_0 - h,
  x_0 + h]) = h^{\beta + 1}$ where $\beta \geq 0$ for any $h$ in a
  neighbourhood of $0$. If
  \begin{equation*}
    I_n = [x_0 - (\ell_n / n)^{1 / (1 + 2s + \beta)}, x_0 + (\ell_n /
    n)^{1 / (1 + 2s + \beta)}],
  \end{equation*}
  we can contruct an estimator $\wh f_n$ such that
  \begin{equation*}
    \sup_{f \in F} \mb E_{f\mu} \big[ w ( \sup_{x \in I_n} r_n(x)^{-1}
    |\wh f_n(x) - f(x)| ) \big] = o(1).
  \end{equation*}
\end{proposition}

This proposition entails that $r_n(\cdot)$ can be improved for
localized risks~\eqref{eq:improve_r_n} over intervals $I_n$ with size
$(\ell_n / n)^{1 / (1 + 2s + \beta)}$ where $\ell_n$ can be a slow
term such has $(\log n)^{\gamma}$ for any $\gamma \geq 0$. A
consequence is that the lower bound in Theorem~\ref{thm:lower_bound}
cannot be improved, since~\eqref{eq:lower_bound} does not hold anymore
when $I_n$ has a length smaller than~\eqref{eq:I_n_length}. This
phenomenon is linked both to the choice of the uniform metric for
measuring the error of estimation, and to the nature of the noise
within the model~\eqref{eq:model}. It is also a consequence of the
minimax paradigm: it is well-known that the minimax risk actually
concentrates on some critical functions of the considered class (that
we rescale and place within $I_n$ here, hence the critical length for
$I_n$), which is a property allowing to prove lower bounds such as the
one in Theorem~\ref{thm:lower_bound}.

\section{Construction of an adaptive estimator}
\label{sec:estimator}

The adaptive method proposed here differs from the techniques
mentioned in Introduction. Indeed, it is not appropriate here to apply
a wavelet decomposition of the scaling coefficients at the finest
scale since it is a $\mbb L^2$-transform, while the
criterion~\eqref{eq:new_UB} considered here uses the uniform
metric. This is the reason why we focus the analysis on a precise
estimation of the scaling coefficients. The technique consists in a
local polynomial approximation of $f$ within adaptively selected
bandwidths for each scaling coefficient.

Let $(V_j)_{j \geq 0}$ be a multiresolution analysis of $\mb L^2([0,
1])$ with scaling function $\phi$ compactly supported and $R$-regular
(the parameter $R$ comes from Theorem~\ref{thm:upper_bound}), which
ensures that
\begin{equation}
  \label{eq:unif_approx}
  \norminfty{f - P_j f} \lesssim 2^{-j s}
\end{equation}
for any $f \in H(s, L)$ with $s \in (0, R+1]$, where $P_j$ denotes the
projection onto $V_j$. We use $P_j$ as an interpolation
transform. Interpolation transforms in the unit interval are
constructed in \cite{donoho92} and \cite{cohen_daubechies_vial93}. We
have $P_j f = \sum_{k = 0}^{2^j - 1} \alpha_{jk} \phi_{jk}$, where
$\phi_{jk}(\cdot) = 2^{j/2} \phi(2^j \cdot - k)$ and $\alpha_{jk} =
\int f \phi_{jk}$. We consider the largest integer $J$ such that $N :=
2^J \leq n$, and we estimate the scaling coefficients at the high
resolution $J$. For appropriate estimators $\wh \alpha_{Jk}$ of
$\alpha_{Jk}$, we simply consider
\begin{equation}
  \label{eq:estimator}
  \wh f_n := \sum_{k=0}^{2^J - 1} \wh \alpha_{Jk} \phi_{Jk}.
\end{equation}
Let us denote by $\text{Pol}_R$ the set of all real polynomials with
degree at most $R$. If $\bar f_k \in \text{Pol}_R$ is close to $f$
over the support of $\phi_{Jk}$, then
\begin{equation*}
  \alpha_{Jk} = \textstyle \int f \phi_{Jk} \approx \int \bar f_k
  \phi_{Jk}.
\end{equation*}
When the scaling function $\phi$ has $R$ moments, that is
\begin{equation}
  \label{eq:moments}
  \int \phi(t) t^p dt = \ind{p = 0}, \quad p \in \{ 0, \ldots, R \},
\end{equation}
and when $f$ is $s$-H\"older for $s \in (0, R+1]$, accurate estimators
of $\wh \alpha_{Jk}$ are given by
\begin{equation}
  \label{eq:coeff_estim}
  \wh \alpha_{Jk} := 2^{-J/2} \bar f_k(k 2^{-J}).
\end{equation}
If $\phi$ does not satisfies~\eqref{eq:moments}, $\int \bar f
\phi_{Jk}$ can be computed exactly using a quadrature formula, in the
same way as in~\cite{delyon_juditsky95}. Indeed, there is a matrix
$\msf Q_J$ (characterized by $\phi$) with entries $(q_{Jkm})$ for $(k,
m) \in \{ 0, \ldots, 2^J-1 \}^2$ such that
\begin{equation}
  \label{eq:quadrature}
  \int P \phi_{Jk} = 2^{-J / 2}\sum_{m \in \Gamma_{Jk}}
  q_{Jkm} P(m / 2^{J})
\end{equation}
for any $P \in \text{Pol}_R$. Within this equation, the entries of the
quadrature matrix $\msf Q_J$ satisfy
\begin{equation}
  \label{eq:Q_matrix_prop}
  q_{Jkm} \neq 0 \rightarrow |k - m| \leq L_{\phi} \text{ and } m \in
  \Gamma_{Jk},
\end{equation}
where $L_{\phi} > 0$ is the support length of $\phi$. Therefore, the
matrix $\msf Q_J$ is band-limited. For instance, if we consider the
Coiflets basis, which satisfies the moment
condition~\eqref{eq:moments}, we have $q_{Jkm} = \ind{k = m}$, and we
can use directly~\eqref{eq:coeff_estim}. If the $(\phi(\cdot - k))_k$
are orthogonal, then $q_{Jkm} = \phi(m - k)$, see
\cite{delyon_juditsky95}.

For the sake of simplicity, we assume in what follows that $\phi$
satisfies the moment condition~\eqref{eq:moments}, thus $\alpha_{Jk}$
is estimated by~\eqref{eq:coeff_estim}. Each polynomial $\bar f_k$
in~\eqref{eq:coeff_estim} is defined via a least square minimization
which is localized within a data-driven bandwidth $\wh \Delta_k$,
hence
\begin{equation*}
  \bar f_k = \bar f_k^{(\wh \Delta_k)}.
\end{equation*}
Below, we describe the computation of these polynomials and then, we
define the selection rule for the $\wh \Delta_k$.

\subsection*{Local polynomials}

The polynomials used to estimate each scaling coefficients are defined
via a slightly modified version of the local polynomial estimator
(LPE). This linear method of estimation is standard, see for instance
\cite{fan_gijbels95, fan_gijbels96}, among many others. For any
interval $\delta \subset [0, 1]$, we define the empirical sample
measure
\begin{equation*}
  \bar \mu_n(\delta) := \frac{1}{n} \sumin \ind{\delta}(X_i),
\end{equation*} 
where $\ind{\delta}$ is the indicator of $\delta$, and if $\bar
\mu_n(\delta) > 0$, we introduce the pseudo-inner product
\begin{equation}
  \label{eq:design_prodsca}
  \prodsca{f}{g}_{\delta} := \frac{1}{\bar \mu_n(\delta)}
  \int_{\delta} f g \, d \bar \mu_n,
\end{equation}
and $\norm{g}_{\delta} := \prodsca{g}{g}_{\delta}^{1/2} $ the
corresponding pseudo-norm. The LPE consists in looking for the
polynomial $\bar f^{(\delta)}$ of degree $R$ which is the closest to
the data in the least square sense, with respect to the localized
design-adapted norm $\norm{\cdot}_{\delta}$:
\begin{equation}
  \label{eq:least_squares}
  \bar f^{(\delta)} := \argmin_{g \in \text{Pol}_R} \norm{Y - g}_{\delta}^2,
\end{equation}
where we recall that $\text{Pol}_R$ is the set of all real polynomials
with degree at most $R$. We can rewrite~\eqref{eq:least_squares} in a
variational form, in which we look for $\bar f^{(\delta)} \in
\text{Pol}_R$ such that for any $\varphi \in \text{Pol}_R$,
\begin{equation}
  \label{eq:variational}
  \prodsca{\bar f^{(\delta)}}{\varphi}_{\delta} =
  \prodsca{Y}{\varphi}_{\delta},
\end{equation}
where it suffices to consider only power functions
$\varphi_{kp}(\cdot) = (\cdot - k / 2^{J})^p$, $0 \leq p \leq R$ when
estimating in a neighbourhood of the regular sampling point $k /
2^{J}$. The coefficients vector $\bar \tta_k^{(\delta)} \in
\setR^{R+1}$ of the polynomial $\bar f_k^{(\delta)}$ is therefore
solution, when it makes sense, of the linear system
\begin{equation*}
  \mb X_k^{(\delta)} \tta = \mb Y_k^{(\delta)},
\end{equation*}
where for $0 \leq p,q \leq R$:
\begin{equation}
  \label{eq:stiff_matrix}
  (\mb X_k^{(\delta)})_{p, q} :=
  \prodsca{\varphi_{kp}}{\varphi_{kq}}_{\delta} \quad \text{ and
  } \quad (\mb Y_k^{(\delta)})_p := \prodsca{Y}{\varphi_{kp}}_{\delta}. 
\end{equation}
We modify this system as follows: when the smallest eigenvalue of $\mb
X_k^{(\delta)}$ (which is non-negative) is too small, we add a
correcting term allowing to bound it from below. We introduce
\begin{equation*}
  \mb {\bar X}_k^{(\delta)} := \mb X_k^{(\delta)} + (n \bar
  \mu_n(\delta))^{-1/2} \mb{Id}_{R + 1}
  \ind{\Omega_k(\delta)^{\complement}},
\end{equation*}
where $\mb{Id}_{R+1}$ is the identity matrix in $\setR^{R+1}$ and
\begin{equation}
 \label{eq:omega_I_def}
 \Omega_k(\delta) := \big\{ \lba(\mb X_k^{(\delta)}) > (n \bar
 \mu_n(\delta))^{-1/2} \big\},
\end{equation}
where $\lba(M)$ stands for the smallest eigenvalue of a matrix
$M$. The quantity $(n \bar \mu_{n}(\delta))^{-1/2}$ comes from the
variance of $\bar f_k^{(\delta)}$, and this particular choice
preserves the convergence rate of the method. This modification of the
classical LPE is convenient in situations with little data.

\begin{definition}
  When $\bar \mu_n(\delta) > 0$, we consider the solution $\bar
  \tta_k^{(\delta)}$ of the linear system
  \begin{equation}
    \label{eq:new_system}
    \mb{\bar X}_k^{(\delta)} \tta = \mb Y_k^{(\delta)},
  \end{equation}
  and introduce $ \bar f_k^{(\delta)}(x) := (\bar \tta_k^{(\delta)})_0
  + (\bar \tta_k^{(\delta)})_1 (x - k / 2^J) + \cdots + (\bar
  \tta_k^{(\delta)})_R (x - k / 2^J)^R$. When $\bar \mu_n(\delta) =
  0$, we take simply $\bar f_k^{(\delta)} := 0$.
\end{definition}

\subsection*{Adaptive bandwidth selection}

The adaptive procedure selecting the intervals $\wh \Delta_k$ is based
on a method introduced by~\cite{lepski90}, see
also~\cite{lepski_mammen_spok_97}, and \cite{lepski_spok97}. If a
family of linear estimators can be ``well-sorted'' by their respective
variances (e.g. kernel estimators in the white noise model, see
\cite{lepski_spok97}), the Lepski procedure selects the largest
bandwidth such that the corresponding estimator does not differ
``significantly'' from estimators with a smaller bandwidth. Following
this principle, we construct a method which adapts to the unknown
smoothness, and additionally to the original Lepski method, to the
distribution of the data (the design density is unknown). Bandwidth
selection procedures in local polynomial estimation can be found in
\cite{fan_gijbels95}, \cite{goldenshluger_nemirovski97} or
\cite{spok98}.

The idea of the adaptive procedure is the following: when $\bar
f^{(\delta)}$ is close to $f$ (that is, when $\delta$ is well-chosen),
we have in view of~\eqref{eq:variational}
\begin{equation*}
  \prodsca{\bar f^{(\delta')} - \bar f^{(\delta)}}{\varphi}_{\delta'}
  = \prodsca{Y - \bar f^{(\delta)}}{\varphi}_{\delta'} \approx
  \prodsca{Y - f}{\varphi}_{\delta'} =
  \prodsca{\xi}{\varphi}_{\delta'}
\end{equation*}
for any $\delta' \subset \delta$, $\varphi \in \text{Pol}_R$, where
the right-hand side is a noise term. Then, in order to ``remove'' this
noise, we select the largest $\delta$ such that this noise term
remains smaller than an appropriate threshold, for any $\delta'
\subset \delta$ and $\varphi = \varphi_{k p}$, $p \in \{ 0, \ldots, R
\}$. The bandwidth $\wh \Delta_k$ is selected in a fixed set of
intervals $G_k$ called \emph{grid} (which is defined below) as
follows:
\begin{equation}
  \label{eq:selection_rule}
  \begin{split}
    \wh \Delta_k := \argmax_{\delta \in G_k} \Big\{ \bar \mu_n(\delta)
    \;|\; \forall \delta' \in &G_k, \delta' \subset \delta,\; \forall
    p \in \{ 0, \ldots, R \}, \\
    &| \prodsca{\bar f_k^{(\delta')} - \bar
      f_k^{(\delta)}}{\varphi_{kp}}_{\delta'} | \leq
    \norm{\varphi_{kp}}_{\delta'} T_n(\delta, \delta') \Big\},
  \end{split}
\end{equation}
where
\begin{equation}
  \label{eq:threshold}
  T_n(\delta, \delta') := \sigma \Big[ \Big( \frac{\log n}{n \bar
    \mu_n(\delta)} \Big)^{1/2} + D C_R \Big( \frac{\log( n \bar
    \mu_n(\delta))}{n \bar \mu_n(\delta') } \Big)^{1/2} \Big],
\end{equation}
with $C_R := 1 + (R + 1)^{1/2}$ and $D > (2(b+1))^{1/2}$, if we want
to prove Theorem~\ref{thm:upper_bound} with a loss function satisfying
$w(x) \lesssim (1 + |x|^b)$. The threshold choice~\eqref{eq:threshold}
can be understood in the following way: since the variance of $\bar
f_k^{(\delta)}$ is of order $(n \bar \mu_n(\delta))^{-1/2}$, we see
that the two terms in $T_n(\delta, \delta')$ are ratios between a
penalizing log term and the variance of the estimators compared by the
rule~\eqref{eq:selection_rule}. The penalization term is linked with
the number of comparisons necessary to select the bandwidth. To prove
Theorem~\ref{thm:upper_bound}, we use the grid
\begin{equation}
  \label{eq:sym_big_grid}
  G_k := \bigcup_{1 \leq i \leq n} \Big\{ \big[ k 2^{-J} - |X_i - k
  2^{-J}|, k 2^{-J} + |X_i - k 2^{—J}| \big] \Big\},
\end{equation}
and we recall that the scaling coefficients are estimated by
\begin{equation*}
  \wh \alpha_{Jk} := 2^{-J/2} \bar f_k^{(\wh \Delta_k)}(k 2^{-J}).
\end{equation*}

\begin{remark}
  In this form, the adaptive estimator has a complexity $O(n^2)$.
  This can be decreased using a smaller grid. An example of such a
  grid is the following: first, we sort the $(X_i, Y_i)$ into
  $(X_{(i)}, Y_{(i)})$ such that $X_{(i)} < X_{(i+1)}$. Then, we
  consider $i(k)$ such that $k / 2^J \in [X_{(i(k))}, X_{(i(k) + 1)}]$
  (if necessary, we take $X_{(0)} = 0$ and $X_{(n+1)} = 1$) and for
  some $a > 1$ (to be chosen by the statistician) we introduce
  \begin{equation}
    \label{eq:geom_grid}
    G_k := \bigcup_{p=0}^{[\log_a (i(k)+1)]} \bigcup_{q = 0}^{[ \log_a
      (n - i(k))]} \Big\{ \big[ X_{(i(k) + 1 - [a^p])} , X_{(i(k) +
      [a^q])}\big] \Big\}.
  \end{equation}
  With this grid, the selection of the bandwidth is fast, and the
  complexity of the procedure is $O(n (\log n)^2)$. We can use this
  grid in practice, but we need extra assumptions on the design if we
  want to prove Theorem~\ref{thm:upper_bound} with this grid choice.
\end{remark}

\section{Proofs}

We recall that the weight function $w(\cdot)$ is non-negative,
non-decreasing and such that $w(x) \leq A (1 + |x|)^b$ for some $A, b
> 0$. We denote by $\mu^n$ the joint law of $X_1, \ldots, X_n$ and
$\mf X_n$ the sigma-field generated by $X_1, \ldots, X_n$. $|A|$
denotes both the length of an interval $A$ and the cardinality of a
finite set $A$. $M^{\top}$ is the transpose of $M$, and $\xi = (\xi_1,
\ldots, \xi_n)^{\top}$.

\subsection*{Proof of Theorem~\ref{thm:upper_bound}}

To prove the upper bound, we use the estimator defined
by~\eqref{eq:estimator} where $\phi$ is a scaling function
satisfying~\eqref{eq:moments} (for instance the Coiflets basis), and
where the scaling coefficients are estimated
by~\eqref{eq:coeff_estim}. Using together~\eqref{eq:unif_approx} and
the fact that $r_n(x) \gtrsim (\log n / n)^{s / (1 + 2s)}$ for any
$x$, we have $\sup_{x \in [0, 1]} r_n(x)^{-1} \norminfty{f - P_Jf} =
o(1)$. Hence,
\begin{align*}
  \sup_{x \in [0, 1]} r_n(x)^{-1} |\wh f_n(x) - f(x)| &\lesssim
  \sup_{x \in [0, 1]} r_n(x)^{-1} \Big| \sum_{k = 0}^{2^J -
    1} (\wh \alpha_{Jk} - \alpha_{Jk}) \phi_{Jk}(x) \Big| \\
  & \lesssim \max_{0 \leq k \leq 2^J - 1} \sup_{x \in S_k} r_n(x)^{-1}
  2^{J/2} |\wh \alpha_{Jk} - \alpha_{Jk}|,
\end{align*}
where $S_k$ denotes the support of $\phi_{Jk}$. Then, expanding $f$ up
to the degree $\ppint{s} \leq R$ and using~\eqref{eq:moments}, we
obtain
\begin{equation}
  \label{eq:unif_control_rate}
  \sup_{x \in [0, 1]} r_n(x)^{-1} |\wh f_n(x) - f(x)| \lesssim \max_{0
    \leq k \leq 2^J - 1} \sup_{x \in S_k} r_n(x)^{-1} |\bar f_k^{(\wh
    \Delta_k)}(x_k) - f(x_k)|.
\end{equation}
Since $|S_k| = 2^{-J} \asymp n^{-1}$, we have
\begin{equation}
  \label{eq:r_n_control}
  \sup_{x \in S_k} r_n(x)^{-1} \lesssim r_n(x_k)^{-1}.
\end{equation}
Indeed, since $\mu$ is continuous, $r_n(\cdot)$ is continuously
differentiable and we have $ \sup_{x \in S_k} | r_n(x)^{-1} -
r_n(x_k)^{-1}| \leq 2^{-J} \norminfty{(r_n^{-1})'}$, where $g'$ stands
for the derivative of $g$. Moreover, $|(r_n(x)^{-1})'| \lesssim
h_n'(x) h_n(x)^{-(s+1)} \lesssim n^{-1}$, since $h_n'(x) \lesssim 1$
and $h_n(x) \gtrsim (\log n / n)^{1 / (2s + 1)}$,
thus~\eqref{eq:r_n_control}.

In what follows, $\norminfty{\cdot}$ denotes the supremum norm in
$\setR^{R+1}$. The following lemma is a version of the bias-variance
decomposition of the local polynomial estimator, which is classical:
see for instance \cite{fan_gijbels95,fan_gijbels96},
\cite{goldenshluger_nemirovski97}, \cite{spok98}, among others. We
define the matrix
\begin{equation*}
  \mb E_k^{(\delta)} := \mb \Lambda_k^{(\delta)} \mb{ \bar
    X}_k^{(\delta)} \mb \Lambda_k^{(\delta)}, 
\end{equation*}
where $\bar {\mb X}_k$ is given by~\eqref{eq:stiff_matrix} and $\mb
\Lambda_k^{(\delta)} := \diag [ \norm{\varphi_{k0}}_{\delta}^{-1},
\ldots, \norm{\varphi_{kR}}_{\delta}^{-1}]$.

\begin{lemma}
  \label{lem:bias_variance}
  Conditionally on $\mf X_n$, for any $f \in H(s, L)$ and $\delta \in
  G_k$, we have
  \begin{equation*}
    | \bar f_k^{(\delta)}(x_k) - f(x_k) | \lesssim \lambda(\mb
    E_k^{(\delta)})^{-1} \big( L |\delta|^s + \sigma (n \bar
    \mu_n(\delta))^{-1/2} \norminfty{\mb U_k^{(\delta)} \xi} \big)
  \end{equation*}
  on $\Omega_k(\delta)$, where $\mb U_k^{(\delta)}$ is a $\mf
  X_n$-measurable matrix of size $(R + 1) \times (n \bar
  \mu_n(\delta))$ satisfying $\mb U_k^{(\delta)} (\mb
  U_k^{(\delta)})^{\top} = \mb{Id}_{R+1}$.
\end{lemma}

Note that within Lemma~\ref{lem:bias_variance}, the bandwidth $\delta$
can change from one point $x_k$ to another. We denote shortly $\mb U_k
:= \mb U_k^{(\delta_k)}$. Let us define $W := \mb U \xi$ where $\mb U
:= (\mb U_0^{\top}, \ldots, \mb U_{2^J}^{\top})^{\top}$. In view of
Lemma~\ref{lem:bias_variance}, $W$ is conditionally on $\mf X_n$ a
centered Gaussian vector such that $\mb E_{f\mu} [ W_k^2 | \mf X_n] =
1$ for any $k \in \{ 0, \ldots, (R+1) 2^J \}$. We introduce $W^N :=
\max_{0 \leq k \leq (R+1) 2^J} |W_k|$ and the event $\mc W_N := \big\{
| W^N - \mb E[W^N | \mf X_n] | \leq L_W (\log n)^{1/2} \big\}$, where
$L_W > 0$. We recall the following classical results about the
supremum of a Gaussian vector (see for instance
in~\cite{ledoux_talagrand91}):
\begin{equation*}
  \mb E_{f\mu} \big[ W^N | \mf X_n \big] \lesssim (\log N)^{1/2}
  \lesssim (\log n)^{1/2},
\end{equation*}
and
\begin{equation}
  \label{eq:gauss_deviation}
  \mb P_{f\mu} \big[ \mc W_N^{\complement} | \mf X_n \big] \lesssim
  \exp( -L_W^2 (\log n) / 2) = n^{-L_W^2 / 2}.
\end{equation}
Let us define the event
\begin{equation*}
  \mrm T_k := \{ \bar \mu_n(\Delta_k) \leq \bar \mu_n(\wh \Delta_k)
  \}
\end{equation*}
and $R_k := \sigma \big( \frac{\log n}{n \bar \mu_n(\Delta_k)}
\big)^{1/2}$ where the intervals $\Delta_k$ are given by
\begin{equation*}
  \Delta_k := \argmax_{\delta \in G_k} \Big\{ \bar \mu_n(\delta) \;|\; L
  |\delta|^s \leq \sigma \Big( \frac{\log n}{n \bar \mu_n(\delta)}
  \Big)^{1/2} \Big\}.
\end{equation*}
There is an event $\mrm S_n \in \mf X_n$ such that $\mu^n[ \mrm
S_n^{\complement}] = o(1)$ faster than any power of $n$, and such that
$R_k \asymp r_n(x_k)$ and $\lambda(\mb E_k^{(\Delta_k)}) \gtrsim 1$,
uniformly for any $k \in \{0, \ldots, 2^J - 1 \}$. This event is
constructed below. We decompose
\begin{equation*}
  |\bar f_k^{(\wh \Delta_k)}(x_k) - f(x_k)| \leq A_k + B_k + C_k +
  D_k,
\end{equation*}
where
\begin{align*}
  A_k &:= |\bar f_k^{(\wh \Delta_k)}(x_k) - f(x_k)| \ind{\mc
    W_N^{\complement} \cup \mrm S_n^{\complement}}, \\
  B_k &:= |\bar f_k^{(\wh \Delta_k)}(x_k) - f(x_k)| \ind{\mrm
    T_k^{\complement} \cap \mc W_N \cap \mrm S_n}, \\
  C_k &:= |\bar f_k^{(\wh \Delta_k)}(x_k) - \bar
  f_k^{(\Delta_k)}(x_k)| \ind{\mrm T_k \cap \mrm S_n}, \\
  D_k &:= |\bar f_k^{(\Delta_k)}(x_k) - f(x_k)| \ind{\mc W_N \cap \mrm
    S_n}.
\end{align*}
%%% STEP A
\emph{Term $A_k$.}~For any $\delta \in G_k$, we have
\begin{equation}
  \label{eq:bad_case}
  | \bar f_k^{(\delta)} (x_k) | \lesssim (n \bar
  \mu_n(\delta))^{1/2} \norminfty{f} (1 + W^N).
\end{equation}
This inequality is proved below. Using~\eqref{eq:bad_case}, we can
bound
\begin{equation*}
  \mb E_{f\mu} \big[ w \big( \max_{0 \leq k \leq 2^J} r_n(x_k)^{-1}
  |\bar f_k^{(\wh \Delta_k)}(x_k)| \big)| \mf X_n \big]
\end{equation*}
by some power of $n$. Using $\norminfty{f} \leq Q$ together with the
fact that $L_W$ can be arbitrarily large in~\eqref{eq:gauss_deviation}
and since $\mb \mu^n[ S_n^{\complement}] = o(1)$ faster than any power
of $n$, we obtain
\begin{equation*}
  \mb E_{f\mu} \big[ w( \max_{0 \leq k \leq 2^J} r_n(x_k)^{-1} A_k)
  \big] = o(1).
\end{equation*}

%%% STEP D
\noindent
\emph{Term $D_k$.}~Using together Lemma~\ref{lem:bias_variance}, the
definition of $\Delta_k$ and the fact that $W^N \lesssim (\log
n)^{1/2}$ on $\mc W_N$, we have
\begin{equation*}
  |\bar f_k^{(\Delta_k)}(x_k) - f(x_k)| \leq \lambda(\mb
  E_k^{(\Delta_k)})^{-1} R_k (1 + (\log n)^{-1/2} W^N) \lesssim
  \lambda(\mb E_k^{(\Delta_k)})^{-1} r_n(x_k)
\end{equation*}
on $\mc W_N \cap \mrm S_n$, thus
\begin{equation*}
  \mb E_{f\mu} \big[ w( \max_{0 \leq k \leq 2^J} r_n(x_k)^{-1} D_k)
  \big] \lesssim 1.
\end{equation*}
\noindent \emph{Term $C_k$.}~We introduce $G_k(\delta) := \{ \delta'
\in G_k | \delta' \subset \delta \}$ and the following events:
\begin{align*}
  \mc T_k(\delta, \delta', p) &:= \big\{ | \prodsca{\bar
    f_k^{(\delta)} - \bar f_k^{(\delta')}}{\varphi_{kp}}_{\delta'}|
  \leq \sigma \norm{\varphi_{kp}}_{\delta'} T_n(\delta, \delta')
  \big\}, \\
  \mc T_k(\delta, \delta') &:= \cap_{0 \leq p \leq R} \mc T_k(\delta,
  \delta'), \\
  \mc T_k(\delta) &:= \cap_{\delta' \in G_k(\delta)} \mc T_k(\delta,
  \delta').
\end{align*}
By the definition~\eqref{eq:selection_rule} of the selection rule, we
have $\mrm T_k \subset \mc T_k(\wh \Delta_k, \Delta_k)$. Let $\delta
\in G_k, \delta' \in G_k(\delta)$. On $\mc T_k(\delta, \delta') \cap
\Omega_k(\delta')$ we have (see below)
\begin{equation}
  \label{eq:est_diff}
  | \bar f_k^{(\delta)}(x_k) - \bar f_k^{(\delta')}(x_k) | \lesssim
  \lambda(\mb E_k^{(\delta')})^{-1} \Big( \frac{\log n}{n \bar
    \mu_n(\delta')} \Big)^{1/2}.
\end{equation}
Thus, using~\eqref{eq:est_diff}, we obtain
\begin{equation*}
  \mb E_{f\mu} \big[ w( \max_{0 \leq k \leq 2^J} r_n(x_k)^{-1} C_k)
  \big] \lesssim 1.
\end{equation*}
\emph{Term $B_k$.}~By the definition~\eqref{eq:selection_rule} of the
selection rule, we have $\mrm T_k^{\complement} \subset \mc
T_k(\Delta_k)^{\complement}$. We need the following lemma.
\begin{lemma}
  \label{lem:proba_reject}
  If $\delta \in G_k$ satisfies
  \begin{equation}
    \label{eq:bias_variance}
    L |\delta|^s \leq \sigma \Big( \frac{\log n}{n \bar
      \mu_n(\delta)} \Big)^{1/2}
  \end{equation}
  and $f \in H(s, L)$, we have
  \begin{equation*}
    \mb P_{f\mu} \big[\mc T_k(\delta)^{\complement} | \mf X_n\big]
    \leq (R+1) ( n \bar \mu_n(\delta))^{1 -D^2 / 2}
  \end{equation*}
  on $\Omega_k(\delta)$, where $D$ is the constant from the
  threshlod~\eqref{eq:threshold}.
\end{lemma}
Using together Lemma~\ref{lem:proba_reject}, $\norminfty{f} \leq Q$
and~\eqref{eq:bad_case}, we obtain
\begin{equation*}
  \mb E_{f\mu} \big[ w\big( \max_{0 \leq k \leq 2^J} R_k^{-1} |\bar
  f_k^{(\wh \Delta_k)}(x_k) - f(x_k)| \ind{\mrm T_k^{\complement} \cap
    \mc W_N} \big) | \mf X_n \big] \lesssim 1,
\end{equation*}
thus
\begin{equation*}
  \mb E_{f\mu} \big[ w( \max_{0 \leq k \leq 2^J} r_n(x_k)^{-1} B_k)
  \big] \lesssim 1,
\end{equation*}
and Theorem~\ref{thm:upper_bound} follows. \hfill $\square$
%%% END OF PROOF

\subsection*{Proof of Lemma~\ref{lem:bias_variance}}

On $\Omega_k(\delta)$, we have $\mb {\bar X}_k^{(\delta)} = \mb
X_k^{\delta}$, and $\mb \lba(\mb X_k^{(\delta)}) > ( n \bar
\mu_n(\delta))^{-1/2} > 0$, thus $\mb X_k^{(\delta)}$ and $\mb
E_k^{(\delta)}$ are invertible. Let $f_k$ be the Taylor polynomial of
$f$ at $x_k$ up to the order $\ppint{s}$ and $\tta_k \in \setR^{R+1}$
be the coefficient vector of $f_k$. Using $f \in H(s, L)$, we obtain
\begin{align*}
  |\bar f_k^{(\delta)}(x_k) - f(x_k)| &\lesssim |\prodsca{(\mb
    \Lba_k^{(\delta)})^{-1} (\bar \tta_k^{(\delta)} -
    \tta_k)}{e_1}| + |\delta|^s \\
  &= |\prodsca{(\mb E_k^{(\delta)})^{-1} \mb \Lba_k^{(\delta)} \mb
    X_k^{(\delta)} (\bar \tta_k^{(\delta)} - \tta_k)}{e_1} | +
  |\delta|^s.
\end{align*}
In view of~\eqref{eq:variational}, we have on $\Omega_k(\delta)$ for
any $p \in \{ 0, \ldots, R \}$:
\begin{align*}
  (\mb X_k^{(\delta)} (\bar \tta_k^{(\delta)} - \tta_k))_p &=
  \prodsca{\bar f_k^{(\delta)} - f_k}{\varphi_{kp}}_{\delta} \\
  &= \prodsca{Y - f_k}{\varphi_{kp}}_{\delta}
\end{align*}
thus, $\mb X_k^{(\delta)} (\bar \tta_k^{(\delta)} - \tta_k) =
B_k^{(\delta)} + V_k^{(\delta)}$ where $(B_k^{(\delta)})_p :=
\prodsca{f - f_k}{\varphi_{kp}}_{\delta}$ and $(V_k^{(\delta)})_p :=
\prodsca{\xi}{\varphi_{kp}}_{\delta}$, which correspond respectively
to bias and variance terms. Since $f \in H(s, L)$ and $\lba(M)^{-1} =
\norm{M^{-1}}$ for any symmetrical and positive matrix $M$, we have
\begin{equation*}
  |\prodsca{(\mb E_k^{(\delta)})^{-1} \mb \Lba_k^{(\delta)}
    B_k^{(\delta)}}{e_1}| \lesssim \lba(\mb E_k^{(\delta)})^{-1} L
  |\delta|^s.
\end{equation*}
Since $(V_k^{(\delta)})_p = (n \bar \mu_n(\delta))^{-1} \mb
D_k^{(\delta)} \xi$ where $\mb D_k^{(\delta)}$ is the $(R+1) \times (n
\bar \mu_n(\delta))$ matrix with entries $(\mb D_k^{(\delta)})_{i, p}
:= (X_i - x_k)^p$, $X_i \in \delta$, we can write
\begin{equation*}
  | \prodsca{(\mb E_k^{(\delta)})^{-1} \mb \Lambda_k^{(\delta)}
    V_k^{(\delta)}}{e_1}_{\delta} | \lesssim \sigma (n \bar
  \mu_n(\delta))^{-1/2} \norm{(\mb E_k^{(\delta)})^{-1/2}}
  \norminfty{\mb U_k^{(\delta)} \xi},
\end{equation*}
where $\mb U_k^{(\delta)} := (n \bar \mu_n(\delta))^{-1/2} (\mb
E_k^{(\delta)})^{-1/2} \mb \Lambda_k^{(\delta)} \mb D_k^{(\delta)}$
satisfies $\mb U_k^{(\delta)} (\mb U_k^{(\delta)})^{\top} =
\mb{Id}_{R+1}$ since $\mb E_k^{(\delta)} = \mb \Lambda_k^{(\delta)}
\mb X_k^{(\delta)} \mb \Lambda_k^{(\delta)}$ and $\mb X_k^{(\delta)} =
(n \bar \mu_n(\delta))^{-1} \mb D_k^{(\delta)} (\mb
D_k^{(\delta)})^{\top}$, thus the lemma. $\hfill \square$

\subsection*{Proof of~\eqref{eq:bad_case}}

If $\bar \mu_n(\delta) = 0$, we have $\bar f_k^{(\delta)} = 0$ by
definition and the result is obvious, thus we assume $\bar
\mu_n(\delta) > 0$. Since $\lba( \mb {\bar X}_k^{(\delta)}) \geq (n
\bar \mu_n(\delta))^{-1/2} > 0$, $\mb {\bar X}_{k}^{(\delta)}$ and
$\mb \Lba_{k}^{(\delta)}$ are invertible and $\mb E_{k}^{(\delta)}$
also is. The proof of~(\ref{eq:bad_case}) is then similar to that of
Lemma~\ref{lem:bias_variance}, where the bias is bounded by
$\norminfty{f}$ and where we use the fact that $\lba( \mb {\bar
  X}_k^{(\delta)}) \geq (n \bar \mu_n(\delta))^{-1/2}$
to control the variance term. $\hfill \square$

\subsection*{Proof of~\eqref{eq:est_diff}}

Let us define $\mb H_k^{(\delta)} := \mb \Lba_k^{(\delta)} \mb
X_k^{(\delta)}$. On $\Omega_k(\delta)$, we have:
\begin{equation*}
  |\bar f_k^{(\delta)}(x_k) - \bar f_k^{(\delta')}(x_k)| = | (\bar
  \tta_k^{(\delta)} - \bar \tta_k^{(\delta')})_0 |
  \lesssim \lba(\mb E_k^{(\delta')})^{-1} \norminfty{\mb
    H_k^{(\delta')} ( \bar \tta_k^{(\delta)} - \bar \tta_k^{(\delta')}
    ) }.
\end{equation*}
Since on $\Omega_k(\delta')$, $( \mb H_k^{(\delta')} (\bar
\tta_k^{(\delta)} - \bar \tta_k^{(\delta')}) )_p = \prodsca{\bar
  f_k^{(\delta)} - \bar f_k^{(\delta')}}{\varphi_{kp}}_{\delta'} /
\norm{\varphi_{kp}}_{\delta'}$, and since $\delta' \subset \delta$, we
obtain~\eqref{eq:est_diff} on $\mc T_k(\delta, \delta')$. $\hfill
\square$

\subsection*{Proof of Lemma~\ref{lem:proba_reject}}

We denote by $\mb P_k^{(\delta)}$ the projection onto
$\Span\{\varphi_{k0}, \ldots, \varphi_{kR}\}$ with respect to the
inner product $\prodsca{\cdot}{\cdot}_{\delta}$. Note that on
$\Omega_k(\delta)$, we have $\bar f_k^{(\delta)} = \mb P_k^{(\delta)}
Y$. Let $\delta \in G_k$ and $\delta' \in G_k(\delta)$. In view
of~\eqref{eq:variational}, we have on $\Omega_k(\delta)$ for any
$\varphi = \varphi_{kp}$, $p \in \{ 0, \ldots, R \}$:
\begin{align*}
  \prodsca{\bar f_k^{(\delta')} - \bar
    f_k^{(\delta)}}{\varphi}_{\delta'} &= \prodsca{Y - \bar
    f_k^{(\delta)}}{\varphi}_{\delta'} \\
  &= \prodsca{f - \mb P_k^{(\delta)} Y}{\varphi}_{\delta'} +
  \prodsca{\xi}{\varphi}_{\delta'} \\
  &= A_k - B_k + C_k,
\end{align*}
where $A_k := \prodsca{f - \mb P_k^{(\delta)} f}{\varphi}_{\delta'}$,
$B_k := \sigma \prodsca{\mb P_k^{(\delta)} \xi}{\varphi}_{\delta'}$
and $C_k := \sigma \prodsca{\xi}{\varphi}_{\delta'}$. If $f_k$ is the
Taylor polynomial of $f$ at $x_k$ up to the order $\ppint{s}$, since
$\delta' \subset \delta$ and $f \in H(s, L)$ we have:
\begin{align*}
  |A_k| \leq \norm{\varphi}_{\delta'} \norm{f - f_k + \mb
    P_k^{(\delta)} (f_k - f)}_{\delta} \leq \norm{\varphi}_{\delta'}
  \norm{f - f_k}_{\delta} \lesssim \norm{\varphi}_{\delta'} L
  |\delta|^s,
\end{align*}
and using~\eqref{eq:bias_variance}, we obtain $|A_k| \lesssim
\norm{\varphi}_{\delta'} \sigma \big( \frac{\log n}{n \bar
  \mu_n(\delta)} \big)^{1/2}$. Since $\mb P_k^{(\delta)}$ is an
orthogonal projection, the variance of $B_k$ is equal to
\begin{align*}
  \sigma^2 \mb E_{f\mu} \big[ \prodsca{\mb P_k^{(\delta)}
    \xi}{\varphi}_{\delta'}^2 | \mf X_n \big] &\leq \sigma^2
  \norm{\varphi}_{\delta'}^2 \mb E_{f\mu} \big[ \norm{\mb
    P_k^{(\delta)} \xi}_{\delta'}^2 | \mf X_n \big] \\
  &= \sigma^2 \norm{\varphi}_{\delta'}^2 \trace(\mb P_k^{(\delta)}) /
  (n \bar \mu_n(\delta')),
\end{align*}
where $\trace(M)$ stands for the trace of a matrix $M$. Since $\mb
P_k^{(\delta)}$ is the projection onto $\text{Pol}_R$, $\trace(\mb
P_k^{(\delta)}) \leq R + 1$, and the variance of $B_k$ is smaller than
$\sigma^2 \norm{\varphi}_{\delta'}^2 (R + 1) / (n \bar
\mu_n(\delta'))$. Then,
\begin{equation}
  \label{eq:B_plus_C_variance}
  \mb E_{f\mu} [ (B+C)^2 | \mf X_n ] \leq \sigma^2
  \norm{\varphi}_{\delta'}^2 C_R^2 / (n \bar \mu_n(\delta')).
\end{equation}
In view of the threshold choice~\eqref{eq:threshold}, we have
\begin{align*}
  \big\{ | \prodsca{\bar f_k^{(\delta)} - \bar
    f_k^{(\delta')}}{\varphi}_{\delta'} | &> \norm{\varphi}_{\delta'}
  T_n(\delta, \delta') \big\} \\
  &\subset \Big\{ \frac{ \norm{\varphi}_{\delta'}^{-1} |B_k + C_k|
  }{\sigma (n \bar \mu_n(\delta'))^{-1/2} C_{R} } > D \big( \log (n
  \bar \mu_n(\delta)) \big)^{1/2} \Big\},
\end{align*}
and using~\eqref{eq:B_plus_C_variance} together with $\mb P[ |N(0, 1)|
> x ] \leq \exp(-x^2 / 2)$ and $|G_k(\delta)| \leq (n \bar
\mu_n(\delta))$, we obtain
\begin{align*}
  \mb P_{f\mu} [ \mc T(\delta)^{\complement} | \mf X_n] &\leq
  \sum_{\delta' \in G_k(\delta)} \sum_{p=0}^R \exp
  \big( -D^2 \log( n \bar \mu_n(\delta) ) / 2 \big) \\
  &\leq (R + 1) (n \bar \mu_n(\delta))^{1 - D^2 / 2},
\end{align*}
which concludes the proof.  $\hfill \square$

\subsection*{Construction of $\mrm S_n$}

We construct an event $\mrm S_n \in \mf X_n$ such that $\mu^n\big[\mrm
S_n^{\complement}\big] = o(1)$ faster than any power of $n$, and such
that on this event, $R_k \asymp r_n(x_k)$ and $\lambda(\mb
E_k^{(\Delta_k)}) \gtrsim 1$ uniformly for any $k \in \{ 0, \ldots,
2^J \}$. We need preliminary approximation results, linked with the
approximation of $\mu$ by $\bar \mu_n$. The following deviation
inequalities use Berstein inequality for the sum of independent random
variables, which is standard. We have
\begin{equation}
  \label{eq:bernstein1}
  \mu^n \Big[ \Big|  \frac{\bar \mu_n(\delta)}{\mu(\delta)} - 1 \Big|
  \Big] \lesssim \exp\big( -\von^2 n \mu(\delta) \big)
\end{equation}
for any interval $\delta \subset [0, 1]$ and $\von \in (0, 1)$. Let us
define the events
\begin{equation*}
  \mrm D_{n, a}^{(\delta)}(x, \von) := \Big\{ \Big|
  \frac{1}{\mu(\delta)} \int_{\delta} \Big( \frac{\cdot - x}{|\delta|}
  \Big)^a d \bar \mu_n - e_a(x, \mu) \Big| \leq \von \Big\}
\end{equation*}
where $e_a(x, \mu) := (1 + (-1)^a)(\beta(x) + 1) / (a + \beta(x) + 1)$
($a$ is a natural integer) where we recall that $\beta(x)$ comes from
assumption~\ref{ass:design} (if $x$ is such that $\mu(x) > 0$ then
$\beta(x) = 0$). Using together Bernstein inequality and the fact that
\begin{equation*}
  \frac{1}{\mu(\delta)} \int_{\delta} \Big( \frac{t - x}{|\delta|}
  \Big)^a \mu(t) dt \raro e_{a}(x, \mu)
\end{equation*}
as $|\delta| \raro 0$, we obtain
\begin{equation}
  \label{eq:bernstein2}
  \mu^n \big[ (\mrm D_{n, a}^{(\delta)}(x, \von))^{\complement} \big]
  \lesssim \exp \big( -\von^2 n \mu(\delta) \big).
\end{equation}
By definition~\eqref{eq:sym_big_grid} of $G_k$, we have $\Delta_k =
[x_k - H_n(x_k), x_k + H_n(x_k)]$ where
\begin{equation}
  \label{eq:H_def}
  H_n(x) := \argmin_{h \in [0, 1]} \Big\{ L h^s \geq \sigma \Big(
  \frac{\log n}{n \bar \mu_n([x - h, x + h])} \Big)^{1/2} \Big\}
\end{equation}
is an approximation of $h_n(x)$ (see~\eqref{eq:r_n_def}). Since $\bar
\mu_n$ is ``close'' to $\mu$, these quantities are close to each other
for any $x$. Indeed, if $\delta_n(x) := [x - h_n(x), x + h_n(x)]$ and
$\Delta_n(x) := [x - H_n(x), x + H_n(x)]$ we have using
together~\eqref{eq:H_def} and~\eqref{eq:r_n_def}:
\begin{equation}
  \label{eq:trick1}
  \big\{ H_n(x) \leq (1 + \von) h_n(x) \big\} = \Big\{ \frac{\bar
    \mu_n[(1+\von) \delta_n(x)] }{\mu[ \delta_n(x)]} \geq (1 -
  \von)^{-2} \Big\}
\end{equation}
for any $\von \in (0, 1)$, where $(1 + \von) \delta_n(x) := [x -
(1+\von) h_n(x), x + (1 + \von) h_n(x)]$. Hence, for each $x = x_k$,
the left hand side event of~\eqref{eq:trick1} has a probability that
can be controlled under assumption~\ref{ass:design}
by~\eqref{eq:bernstein1}, and the same argument holds for $\{ H_n(x) >
(1 - \von) h_n(x) \}$. Combining~\eqref{eq:bernstein1},
\eqref{eq:bernstein2} and \eqref{eq:trick1}, we obtain that the event
\begin{equation*}
  \mrm B_{n, a}(x, \von) := \Big\{ \Big|
  \frac{1}{\bar \mu_n(\Delta_n(x))} \int_{\Delta_n(x)} \Big(
  \frac{\cdot - x}{|\delta_n(x)|} \Big)^a d \bar \mu_n - e_a(x, \mu)
  \Big| \leq \von \Big\}
\end{equation*}
satisfies also~\eqref{eq:bernstein2} for $n$ large enough. This proves
that $(\mb X_k^{(\Delta_k)})_{p, q}$ and $(\mb
\Lambda_k^{(\Delta_k)})_p$ are close to $e_{p+q}(x_k, \mu)$ and
$e_{2p}(x_k, \mu)^{-1/2}$ respectively on the event
\begin{equation*}
  \mrm S_n := \bigcap_{a \in \{0, \ldots, 2R\}} \bigcap_{k \in \{ 0,
    \ldots, 2^J-1 \} } \mrm B_{n, a}(x_k, \von).
\end{equation*}
Using the fact that $\lambda(M) = \inf_{\norm{x} = 1} x^{\top} M x$
for a symmetrical matrix $M$, where $\lambda(M)$ denotes the smallest
eigenvalue of $M$, we can conclude that for $n$ large enough,
\begin{equation*}
  \lambda(\mb \Lambda_k^{(\Delta_k)} \mb X_k^{(\Delta_k)} \mb
  \Lambda_k^{(\Delta_k)}) \gtrsim \min_{x \in [0, 1]} \lambda(\mb E(x,
  \mu)),
\end{equation*}
where $\mb E(x, \mu)$ has entries $(\mb E(x, \mu))_{p,q} =
e_{p+q}(x,\mu) / (e_{2p}(x, \mu) e_{2q}(x, \mu))^{1/2}$. Since $\mb
E(x,\mu)$ is definite positive for any $x \in [0, 1]$, we obtain that
on $\mrm S_n$, $\lambda(\mb X_k^{(\Delta_k)}) \gtrsim 1$, thus $\mrm
S_n \subset \Omega_n(\Delta_k)$ and $\lambda(\mb E_k^{(\Delta_k)})
\gtrsim 1$ uniformly for any $k \in \{ 0, \ldots, 2^J-1\}$, since $\mb
E_k^{(\Delta_k)} = \mb \Lambda_k^{(\Delta_k)} \mb X_k^{(\Delta_k)} \mb
\Lambda_k^{(\Delta_k)}$ on $\Omega_n(\Delta_k)$. Moreover, since $R_k
= L H_n(x_k)^s$, using together~\eqref{eq:bernstein1} and
\eqref{eq:trick1}, we obtain $R_k \asymp r_n(x_k)$ uniformly for $k\in
\{ 0, \ldots, 2^J - 1\}$. \hfill $\square$

\subsection*{Proof of Theorem~\ref{thm:lower_bound}}

The main features of the proof are first, a reduction to the Bayesian
risk over an hardest cubical subfamily of functions for the $\mbb
L^{\infty}$ metrics, which is standard: see \cite{korostelev93},
\cite{donoho94}, \cite{korostelev_nussbaum_99} and \cite{bertin02},
and the choice of rescaled hypothesis with design-adapted bandwidth
$h_n(\cdot)$, necessary to achieve the rate $r_n(\cdot)$.

Let us consider $\varphi \in H(s, L; \setR)$ (the extension of $H(s,
L)$ to the whole real line) with support $[-1, 1]$ and such that
$\varphi(0) > 0$. We define
\begin{equation*}
  a := \min \Big[1, \Big( \frac{2}{\norminfty{\varphi}^2} \Big(
  \frac{1}{1 + 2s + \beta} - \alpha \Big) \Big)^{1/(2s)} \Big]
\end{equation*}
and
\begin{equation*}
  \Xi_n := 2 a (1 + 2^{1 / (s - \ppint{s})}) \sup_{x \in [0, 1]}
  h_n(x),
\end{equation*}
where we recall that $\ppint{s}$ is the largest integer smaller than
$s$. Note that~\eqref{eq:assump_low_bound} entails
\begin{equation}
  \label{eq:proof_lb_sup_h_n}
  \Xi_n \lesssim (\log n / n)^{1 /  (1 + 2s + \beta)}.
\end{equation}
If $I_n = [c_n, d_n]$, we introduce $ x_k := c_n + k \, \Xi_n$ for $k
\in K_n := \big\{ 1, \ldots, \big[|I_n| \, \Xi_n^{-1} \big] \big\}$,
and denote for the sake of simplicity $h_k := h_n(x_k)$. We consider
the family of functions
\begin{equation*}
  f(\cdot; \theta) := \sum_{k \in K_n} \theta_k f_k(\cdot), \quad
  f_k(\cdot) := L a^s h_k^s \varphi \Big( \frac{\cdot - x_k}{h_k} \Big),
\end{equation*}
which belongs to $H(s,L)$ for any $\theta \in [-1, 1]^{|K_n|}$. Using
Bernstein inequality, we can see that
\begin{equation*}
  \mrm H_n := \bigcap_{k \in K_n} \Big\{ \frac{\bar \mu_n([x_k - h_k,
    x_k + h_k])}{\mu([x_k - h_k, x_k + h_k])} \geq  1/ 2 \Big\}
\end{equation*}
satisfies
\begin{equation}
  \label{eq:prob_H_n}
  \mu^n[ \mrm H_n ] = 1 - o(1).
\end{equation}
Let us introduce $b := c^s \varphi(0)$. For any distribution $\mb B$
on $\Theta_n \subset [-1, 1]^{|K_n|}$, by a minoration of the minimax
risk by the Bayesian risk, and since $w$ is non-decreasing, the left
hand side of~\eqref{eq:lower_bound} is smaller than
\begin{align*}
  w(b) \inf_{\wh \theta} \int_{\Theta_n} &\mb P_{\theta}^n \big[
  \max_{k \in K_n} |\wh \theta_k - \theta_k | \geq 1 \big] \mb
  B(d\theta) \\
  &\geq w(b) \int_{\mrm H_n} \inf_{\wh \theta} \int_{\Theta_n} \mb
  P_{\theta}^n \big[ \max_{k \in K_n} |\wh \theta_k - \theta_k | \geq
  1 | \mf X_n\big] \mb B(d\theta) d \mu^n.
\end{align*}
Hence, together with~\eqref{eq:prob_H_n},
Theorem~\ref{thm:lower_bound} follows if we show that on $\mrm H_n$
\begin{equation}
  \label{eq:proof_lb_1}
  \sup_{\wh \theta} \int_{\Theta_n} \mb P_{\theta}^n \big[ \max_{k \in
    K_n} |\wh \theta_k - \theta_k | < 1 | \mf X_n \big] \mb B(d\theta)
  = o(1).
\end{equation}
We denote by $L(\theta; Y_1, \ldots, Y_n)$ the conditional on $\mf
X_n$ likelihood function of the observations $Y_i$
from~\eqref{eq:model} when $f(\cdot) = f(\cdot;\,
\theta)$. Conditionally on $\mf X_n$, we have
\begin{equation*}
  L(\theta; Y_1, \ldots, Y_n) = \prod_{1 \leq i \leq n}
  g_{\sigma}(Y_i) \prod_{k \in K_n} \frac{g_{v_k}(y_k -
    \theta_k)}{g_{v_k}(y_k)},
\end{equation*}
where $g_{v}$ is the density of $N(0, v^2)$, $v_k^2 := \mb E\{ y_k^2 |
\mf X_n\}$ and
\begin{equation*}
  y_k := \frac{\sumin Y_i f_k(X_i)}{\sumin f_k^2(X_i)}.
\end{equation*}
Thus, choosing
\begin{equation*}
  \mb B := \bigotimes_{k \in K_n} \mb b, \quad \mb b := (\delta_{-1} +
  \delta_{1}) /2, \quad \Theta_n := \{-1, 1\}^{|K_n|},
\end{equation*}
the left hand side of~\eqref{eq:proof_lb_1} is smaller than
\begin{align*}
  \int \frac{\prod_{1 \leq i \leq n} g_{\sigma}(Y_i)}{\prod_{k \in
      K_n} g_{v_k} (y_k)} \Big( \prod_{k \in K_n} \sup_{\wh \theta_k }
  \int_{\{-1, 1\}} \ind{|\wh \theta_k - \theta_k | < 1}\; g_{v_k}(y_k
  - \theta_k) \mb b( d\theta_k) \Big) d Y_1 \times \cdots \times d
  Y_n,
\end{align*}
and $\wh \theta_k = \ind{y_k \geq 0} -\ind{y_k < 0}$ are strategies
reaching the supremum. Then, in~\eqref{eq:proof_lb_1}, it suffices to
take the supremum over estimators $\wh \theta$ with coordinates $\wh
\theta_k \in \{ -1, 1 \}$ measurable with respect to $y_k$ only. Since
conditionally on $\mf X_n$, $y_k$ is in law $N(\theta_k, v_k^2)$, the
left hand side of~\eqref{eq:proof_lb_1} is smaller than
\begin{equation*}
  \prod_{k \in K_n} \Big( 1 - \inf_{\wh \theta_k \in \{-1, 1 \} }
  \int_{\{-1, 1\}} \int \ind{|\wh \theta_k(u) - \theta_k| \geq 1}
  g_{v_k}(u - \theta_k) du \, \mb b(d\theta_k) \Big).
\end{equation*}
Moreover, if $\Phi(x) := \int_{-\infty}^x g_1(t) dt$
\begin{align*}
  \inf_{\wh \theta_k \in \{-1, 1 \} } \int_{\{-1, 1\}} \int &\ind{|\wh
    \theta_k(u) - \theta_k| \geq 1} g_{v_k}(u - \theta_k) du \,
  \mb b(d\theta_k) \\
  &\geq \frac{1}{2} \int \min\big( g_{v_k}(u-1), g_{v_k}(u+1) \big) du
  = \Phi(-1 / v_k).
\end{align*}
On $\mrm H_n$, we have in view of~\eqref{eq:r_n_def}
\begin{align*}
  v_k^2 = \frac{\sigma^2}{\sumin f_k^2(X_i)} \geq \frac{2}{(1 -
    \delta) \norminfty{\varphi}^2 c^{2s} \log n},
\end{align*}
and since $\Phi(-x) \geq \exp(-x^2 / 2) (x \sqrt{2\pi})$ for any $x >
0$, we obtain
\begin{equation*}
  \Phi(-1 / v_k) \gtrsim (\log n)^{-1/2} n^{\{\alpha - 1 / (1 + 2s +
    \beta)\} / 2} =: L_n.
\end{equation*}
Thus, the left hand side of~\eqref{eq:proof_lb_1} is smaller than $(1
- L_n)^{|K_n|}$, and since
\begin{equation*}
  |I_n| \Xi_n^{-1} L_n \gtrsim n^{\{1 / (1 + 2s + \beta) - \alpha\} /
    2} (\log n)^{1/2 - 1/(1 + 2s + \beta)} \raro +\infty
\end{equation*}
as $n \raro +\infty$, Theorem~\ref{thm:lower_bound} follows. \hfill
$\square$

\subsection*{Proof of Corollary~\ref{cor:lower_bound}}

Let us consider the loss function $w(\cdot) = |\cdot|$, and let $\wh
f_n^v$ be an estimator converging with rate $v_n(\cdot)$ over $F$ in
the sense of~\eqref{eq:upper_bound}. Hence,
\begin{align*}
  1 &\lesssim \sup_{f \in F} \mb E_{f\mu} \big[ \sup_{x \in I_n}
  r_n(x)^{-1} |\wh f_n^v(x) - f(x)| \big] \\
  &\leq \sup_{x \in I_n} \frac{v_n(x)}{r_n(x)} \sup_{f \in F} \mb
  E_{f\mu} \big[ \sup_{x \in I_n} v_n(x)^{-1} |\wh f_n^v(x) - f(x)|
  \big] \lesssim \sup_{x \in I_n} \frac{v_n(x)}{r_n(x)},
\end{align*}
where we used Theorem~\ref{thm:lower_bound}. \hfill $\square$

\subsection*{Proof of Proposition~\ref{prop:suroptimal}}

Without loss of generality, we consider the loss $w(\cdot) =
|\cdot|$. For proving Proposition~\ref{prop:suroptimal}, we use the
linear LPE. If we denote by $\partial^m f$ the $m$-th derivative of
$f$, a slight modification of the proof of
Lemma~\ref{lem:bias_variance} gives for $f\in H(s, L)$ with $s > m$,
\begin{equation*}
  | \partial^m \bar f_k^{(\delta)}(x_k) - \partial^m f(x_k) | \lesssim
  \lambda(\mb E_k^{(\delta)})^{-1} |\delta|^{-m} \big( L |\delta|^s +
  \sigma (n \bar \mu_n(\delta))^{-1/2} W^N \big),
\end{equation*}
where in the same way as in the proof of
Theorem~\ref{thm:upper_bound}, $W^N$ satisfies
\begin{equation}
  \label{eq:E_WN}
  \mb E_{f\mu} [W^N | \mf X_n] \lesssim (\log N)^{1/2},
\end{equation}
with $N$ depending on the size of the supremum, to be specified
below. First, we prove a). Since $|I_n| \sim (\ell_n / n)^{1 / (2s +
  1)}$, if $I_n = [a_n, b_n]$, the points
\begin{equation*}
  x_k := a_n + (k / n)^{1 / (2s + 1)}, \quad k \in \{ 0, \ldots, N \},
\end{equation*}
where $N := [\ell_n]$ belongs to $I_n$. We consider the bandwidth
\begin{equation}
  \label{eq:h_def_prop}
  h_n = \Big( \frac{\log \ell_n}{n} \Big)^{1/(2s + 1)},
\end{equation}
and we take $\delta_k := [x_k - h_n, x_k + h_n]$. Note that since
$\mu(x) > 0$ for any $x$, $\bar \mu_n(\delta) \asymp |\delta|$ as
$|\delta| \raro 0$ with probability going to $1$ faster than any power
of $n$ (using Berstein inequality, for instance). We consider the
estimator defined by
\begin{equation}
  \label{eq:prop_est}
  \wh f_n(x) := \sum_{m=0}^{r} \partial^m \bar f_k^{(\delta_k)} (x_k) (x
  - x_k)^m / m! \quad \text{ for } x \in [x_k, x_{k+1}), \quad k \in
  \{ 0, \ldots, [\ell_n] \},
\end{equation}
where $r := \ppint{s}$. Using a Taylor expansion of $f$ up to the
degree $r$ together with~\eqref{eq:h_def_prop} gives
\begin{equation*}
  (n / \log n)^{s / (1 + 2s)} \sup_{x \in I_n} |\wh f_n(x) - f(x)|
  \lesssim  \Big( \frac{\log \ell_n}{\log n} \Big)^{s / (1 + 2s)} (1 +
  (\log \ell_n)^{-1/2} W^N).
\end{equation*}
Then, integrating with respect to $\mb P_{f\mu}(\cdot | \mf X_n)$ and
using~\eqref{eq:E_WN} where $N = [\ell_n]$ entails a), since $\log
\ell_n = o(\log n)$.

The proof of b) is similar to that of a). In this setting, the rate
$r_n(\cdot)$ (see~\eqref{eq:r_n_def}) can be written as $r_n(x) =
(\log n / n)^{\alpha_n(x)}$ for $x$ in $I_n$ (for $n$ large enough)
where $\alpha_n(x_0) = s / (1 + 2s + \beta)$ and $\alpha_n(x) > s / (1
+ 2s + \beta)$ for $x \in I_n - \{ x_0 \}$. We define
\begin{equation*}
  x_{k+1} =
  \begin{cases}
    x_k + n^{-\alpha_n(x_k) / s} &\text{ for } k \in \{ -N, \ldots, -1
    \} \\
    x_k + n^{-\alpha_n(x_{k+1}) / s} &\text{ for } k \in \{ 0,
    \ldots, N \}, \\
  \end{cases}
\end{equation*}
where $N := [\ell_n]$. All the points fit in $I_n$, since $|x_{-N} -
x_N| \leq \sum_{-N \leq k \leq N} n^{-\min(\alpha_n(x_k),
  \alpha_n(x_{k+1}))/s} \leq 2 (\ell_n / n)^{1 / (1 + 2s +
  \beta)}$. We consider the bandwidths
\begin{equation*}
  h_k := (\log \ell_n / n)^{\alpha_n(x_k) / s},
\end{equation*}
and the intervals $\delta_k = [x_k - h_k, x_k + h_k]$. We keep the
same definition~\eqref{eq:prop_est} for $\wh f_n$. Since $x_0$ is a
local extremum of $r_n(\cdot)$, we have in the same way as in the
proof of a) that
\begin{align*}
  \sup_{x \in I_n} r_n(x)^{-1} |\wh f_n(x) - f(x)| \lesssim \Big[
  &\max_{-N \leq k \leq -1} \Big( \frac{\log \ell_n}{\log n}
  \Big)^{\alpha_n(x_{k})} \\
  &+ \max_{0 \leq k \leq N-1} \Big( \frac{\log \ell_n}{\log n}
  \Big)^{\alpha_n(x_{k + 1})} \Big] (1 + (\log \ell_n)^{-1/2} W^N),
\end{align*}
hence
\begin{equation*}
  \mb E_{f\mu} \big[ \sup_{x \in I_n} r_n(x)^{-1} |\wh f_n(x) - f(x)|
  \big] \lesssim \Big( \frac{\log \ell_n}{\log n}
  \Big)^{s / (1 + 2s + \beta)} = o(1),
\end{equation*}
which concludes the proof of Proposition~\ref{prop:suroptimal}. \hfill
$\square$

%% BIBLIOGRAPHY

\small

\bibliographystyle{ims}%
\bibliography{biblio}

\begin{thebibliography}{38}
\expandafter\ifx\csname natexlab\endcsname\relax\def\natexlab#1{#1}\fi
\expandafter\ifx\csname url\endcsname\relax
  \def\url#1{\texttt{#1}}\fi
\expandafter\ifx\csname urlprefix\endcsname\relax\def\urlprefix{URL }\fi
\providecommand{\eprint}[2][]{\url{#2}}

\bibitem[{Almansa et~al.(2003)Almansa, Rouge and
  Jaffard}]{almansa_rouge_jaffard03}
\textsc{Almansa, A.}, \textsc{Rouge, B.} and \textsc{Jaffard, S.} (2003).
\newblock Irregular sampling in satellite images and reconstruction algorithms.
\newblock In \textit{CANUM 2003}. CANUM 2003,
  http://www.math.univ-montp2.fr/canum03/communications/ms/andres.almansa.pdf.

\bibitem[{Antoniadis and Fan(2001)}]{antoniadis_fan01}
\textsc{Antoniadis, A.} and \textsc{Fan, J.~Q.} (2001).
\newblock Regularization of wavelet approximations.
\newblock \textit{Journal of the American Statistical Association}, \textbf{96}
  939--967.

\bibitem[{Antoniadis et~al.(1997)Antoniadis, Gregoire and
  Vial}]{antoniadis_et_al97}
\textsc{Antoniadis, A.}, \textsc{Gregoire, G.} and \textsc{Vial, P.} (1997).
\newblock Random design wavelet curve smoothing.
\newblock \textit{Statistics and Probability Letters}, \textbf{35} 225--232.

\bibitem[{Antoniadis and Pham(1998)}]{antoniadis_et_al98}
\textsc{Antoniadis, A.} and \textsc{Pham, D.~T.} (1998).
\newblock Wavelet regression for random or irregular design.
\newblock \textit{Comput. Statist. Data Anal.}, \textbf{28} 353--369.

\bibitem[{Baraud(2002)}]{baraud02}
\textsc{Baraud, Y.} (2002).
\newblock Model selection for regression on a random design.
\newblock \textit{ESAIM Probab. Statist.}, \textbf{6} 127--146 (electronic).

\bibitem[{Bertin(2004)}]{bertin02}
\textsc{Bertin, K.} (2004).
\newblock Minimax exact constant in sup-norm for nonparametric regression with
  random design.
\newblock \textit{J. Statist. Plann. Inference}, \textbf{123} 225--242.

\bibitem[{Brown and Cai(1998)}]{cai_brown98}
\textsc{Brown, L.} and \textsc{Cai, T.} (1998).
\newblock Wavelet shrinkage for nonequispaced samples.
\newblock \textit{The Annals of Statistics}, \textbf{26} 1783--1799.

\bibitem[{Cai and Low(2005)}]{cai_low05}
\textsc{Cai, T.~T.} and \textsc{Low, M.~G.} (2005).
\newblock Nonparametric estimation over shrinking neighborhoods:
  superefficiency and adaptation.
\newblock \textit{Ann. Statist.}, \textbf{33} 184--213.

\bibitem[{Cohen et~al.(1993)Cohen, Daubechies and
  Vial}]{cohen_daubechies_vial93}
\textsc{Cohen, A.}, \textsc{Daubechies, I.} and \textsc{Vial, P.} (1993).
\newblock Wavelets on the interval and fast wavelets transforms.
\newblock \textit{Appl. Comput. Harmon. Anal.}, \textbf{1} 54--81.

\bibitem[{Delouille(2002)}]{delouille_phd}
\textsc{Delouille, V.} (2002).
\newblock \textit{Nonparametric stochastic regression using design-adapted
  wavelets}.
\newblock Ph.D. thesis, Universit\'e catholique de Louvain.

\bibitem[{Delouille et~al.(2001)Delouille, Franke and von
  Sachs}]{delouille_et_al01}
\textsc{Delouille, V.}, \textsc{Franke, J.} and \textsc{von Sachs, R.} (2001).
\newblock Nonparametric stochastic regression with design-adapted wavelets.
\newblock \textit{Sankhy\=a Ser. A}, \textbf{63} 328--366.
\newblock Special issue on wavelets.

\bibitem[{Delouille et~al.(2004)Delouille, Simoens and
  Von~Sachs}]{delouille_et_al04}
\textsc{Delouille, V.}, \textsc{Simoens, J.} and \textsc{Von~Sachs, R.} (2004).
\newblock Smooth design-adapted wavelets for nonparametric stochastic
  regression.
\newblock \textit{Journal of the American Statistical Society}, \textbf{99}
  643--658.

\bibitem[{Delyon and Juditsky(1995)}]{delyon_juditsky95}
\textsc{Delyon, B.} and \textsc{Juditsky, A.} (1995).
\newblock Estimating wavelet coefficients.
\newblock In \textit{Lecture notes in Statistics} (A.~Antoniadis and
  G.~Oppenheim, eds.), vol. 103. Springer-Verlag, New York, 151--168.

\bibitem[{Donoho(1992)}]{donoho92}
\textsc{Donoho, D.} (1992).
\newblock Interpolating wavelet tranforms.
\newblock Tech. rep., Department of Statistics, Stanford University,
  http://www-stat.stanford.edu/~donoho/Reports/1992/interpol.ps.Z.

\bibitem[{Donoho(1994)}]{donoho94}
\textsc{Donoho, D.~L.} (1994).
\newblock Asymptotic minimax risk for sup-norm loss: Solution via optimal
  recovery.
\newblock \textit{Probability Theory and Related Fields}, \textbf{99} 145--170.

\bibitem[{Fan and Gijbels(1995)}]{fan_gijbels95}
\textsc{Fan, J.} and \textsc{Gijbels, I.} (1995).
\newblock Data-driven bandwidth selection in local polynomial fitting: variable
  bandwidth and spatial adaptation.
\newblock \textit{Journal of the Royal Statistical Society. Series B.
  Methodological}, \textbf{57} 371--394.

\bibitem[{Fan and Gijbels(1996)}]{fan_gijbels96}
\textsc{Fan, J.} and \textsc{Gijbels, I.} (1996).
\newblock \textit{Local polynomial modelling and its applications}.
\newblock Monographs on Statistics and Applied Probability, Chapman \& Hall,
  London.

\bibitem[{Feichtinger and Gr{\"o}chenig(1994)}]{grochenig94}
\textsc{Feichtinger, H.~G.} and \textsc{Gr{\"o}chenig, K.} (1994).
\newblock Theory and practice of irregular sampling.
\newblock In \textit{Wavelets: mathematics and applications}. Stud. Adv. Math.,
  CRC, Boca Raton, FL, 305--363.

\bibitem[{Ga\"iffas(2005{\natexlab{a}})}]{gaiffas04a}
\textsc{Ga\"iffas, S.} (2005{\natexlab{a}}).
\newblock Convergence rates for pointwise curve estimation with a degenerate
  design.
\newblock \textit{Mathematical Methods of Statistics}, \textbf{1} 1--27.
\newblock Available at \href{http://hal.ccsd.cnrs.fr/ccsd-00003086/en/}{
  http://hal.ccsd.cnrs.fr/ccsd-00003086/en/ }.

\bibitem[{Ga\"iffas(2005{\natexlab{b}})}]{gaiffas05}
\textsc{Ga\"iffas, S.} (2005{\natexlab{b}}).
\newblock On pointwise adaptive curve estimation based on inhomogeneous data.
\newblock Preprint LPMA no 974 available at
  \href{http://hal.ccsd.cnrs.fr/ccsd-00004605/en/}{
  http://hal.ccsd.cnrs.fr/ccsd-00004605/en/}.

\bibitem[{Goldenshluger and Nemirovski(1997)}]{goldenshluger_nemirovski97}
\textsc{Goldenshluger, A.} and \textsc{Nemirovski, A.} (1997).
\newblock On spatially adaptive estimation of nonparametric regression.
\newblock \textit{Mathematical Methods of Statistics}, \textbf{6} 135--170.

\bibitem[{Guerre(1999)}]{guerre99}
\textsc{Guerre, E.} (1999).
\newblock Efficient random rates for nonparametric regression under arbitrary
  designs.
\newblock Personal communication.

\bibitem[{Hall et~al.(1997)Hall, Marron, Neumann and
  Tetterington}]{hall_et_al97}
\textsc{Hall, P.}, \textsc{Marron, J.~S.}, \textsc{Neumann, M.~H.} and
  \textsc{Tetterington, D.~M.} (1997).
\newblock Curve estimation when the design density is low.
\newblock \textit{The Annals of Statistics}, \textbf{25} 756--770.

\bibitem[{Hall et~al.(1998)Hall, Park and Turlach}]{hall_park_turlach98}
\textsc{Hall, P.}, \textsc{Park, B.~U.} and \textsc{Turlach, B.~A.} (1998).
\newblock A note on design transformation and binning in nonparametric curve
  estimation.
\newblock \textit{Biometrika}, \textbf{85} 469--476.

\bibitem[{Jansen et~al.(2004)Jansen, Nason and Silverman}]{jansen04}
\textsc{Jansen, M.}, \textsc{Nason, P.~G.} and \textsc{Silverman, B.~W.}
  (2004).
\newblock Multivariate nonparametric regression using lifting.
\newblock Tech. rep., University of Bristol, UK,
  http://www.stats.ox.ac.uk/\~{}silverma/pdf/jansennasonsilverman.pdf.

\bibitem[{Kerkyacharian and Picard(2004)}]{kerk_picard_warped_04}
\textsc{Kerkyacharian, G.} and \textsc{Picard, D.} (2004).
\newblock Regression in random design and warped wavelets.
\newblock \textit{Bernoulli}, \textbf{10} 1053--1105.

\bibitem[{Korostelev and Nussbaum(1999)}]{korostelev_nussbaum_99}
\textsc{Korostelev, A.} and \textsc{Nussbaum, M.} (1999).
\newblock The asymptotic minimax constant for sup-norm loss in nonparametric
  density estimation.
\newblock \textit{Bernoulli}, \textbf{5} 1099--1118.

\bibitem[{Korostelev(1993)}]{korostelev93}
\textsc{Korostelev, V.} (1993).
\newblock An asymptotically minimax regression estimator in the uniform norm up
  to exact contant.
\newblock \textit{Theory of Probability and its Applications}, \textbf{38}
  737--743.

\bibitem[{Ledoux and Talagrand(1991)}]{ledoux_talagrand91}
\textsc{Ledoux, M.} and \textsc{Talagrand, M.} (1991).
\newblock \textit{Probability in {B}anach spaces}, vol.~23 of
  \textit{Ergebnisse der Mathematik und ihrer Grenzgebiete (3) [Results in
  Mathematics and Related Areas (3)]}.
\newblock Springer-Verlag, Berlin.
\newblock Isoperimetry and processes.

\bibitem[{Lepski(1990)}]{lepski90}
\textsc{Lepski, O.~V.} (1990).
\newblock On a problem of adaptive estimation in {G}aussian white noise.
\newblock \textit{Theory of Probability and its Applications}, \textbf{35}
  454--466.

\bibitem[{Lepski et~al.(1997)Lepski, Mammen and
  Spokoiny}]{lepski_mammen_spok_97}
\textsc{Lepski, O.~V.}, \textsc{Mammen, E.} and \textsc{Spokoiny, V.~G.}
  (1997).
\newblock Optimal spatial adaptation to inhomogeneous smoothness: an approach
  based on kernel estimates with variable bandwidth selectors.
\newblock \textit{The Annals of Statistics}, \textbf{25} 929--947.

\bibitem[{Lepski and Spokoiny(1997)}]{lepski_spok97}
\textsc{Lepski, O.~V.} and \textsc{Spokoiny, V.~G.} (1997).
\newblock Optimal pointwise adaptive methods in nonparametric estimation.
\newblock \textit{The Annals of Statistics}, \textbf{25} 2512--2546.

\bibitem[{Maxim(2003)}]{voichitaphd}
\textsc{Maxim, V.} (2003).
\newblock \textit{Restauration de signaux bruit\'es sur des plans d'experience
  al\'eatoires}.
\newblock Ph.D. thesis, Universit\'e Joseph Fourier, Grenoble 1.

\bibitem[{Pensky and Wiens(2001)}]{pensky_wiens01}
\textsc{Pensky, M.} and \textsc{Wiens, D.~P.} (2001).
\newblock On non-equally spaced wavelet regression.
\newblock \textit{Advances in Soviet Mathematics}, \textbf{53} 681--690.

\bibitem[{Spokoiny(1998)}]{spok98}
\textsc{Spokoiny, V.~G.} (1998).
\newblock Estimation of a function with discontinuities via local polynomial
  fit with an adaptive window choice.
\newblock \textit{The Annals of Statistics}, \textbf{26} 1356--1378.

\bibitem[{Stone(1982)}]{stone82}
\textsc{Stone, C.~J.} (1982).
\newblock Optimal global rates of convergence for nonparametric regression.
\newblock \textit{The Annals of Statistics}, \textbf{10} 1040--1053.

\bibitem[{V\`azquez et~al.(2000)V\`azquez, Konrad and Dubois}]{konrad00}
\textsc{V\`azquez, C.}, \textsc{Konrad, J.} and \textsc{Dubois, E.} (2000).
\newblock Wavelet-based reconstruction of irregularly-sampled images:
  application to stereo imaging.
\newblock In \textit{Proc. Int. Conf. on Image Processing, ICIP-2000}. IEEE,
  http://iss.bu.edu/jkonrad/Publications/local/cpapers/Vazq00icip.pdf.

\bibitem[{Wong and Zheng(2002)}]{wong_zheng02}
\textsc{Wong, M.-Y.} and \textsc{Zheng, Z.} (2002).
\newblock Wavelet threshold estimation of a regression function with random
  design.
\newblock \textit{Journal of Multivariate Analysis}, \textbf{80} 256--284.

\end{thebibliography}

\end{document}